\documentclass[ejs,preprint]{imsart}

\usepackage{bm}
\usepackage{CJK}
\usepackage{mathrsfs}
\usepackage{mathtools}
\usepackage{extarrows}

\usepackage{leftidx}
\usepackage[utf8]{inputenc}
\usepackage{float}
\usepackage[english]{babel}

\usepackage{amsthm,amsmath}
\RequirePackage[numbers]{natbib}
\RequirePackage[colorlinks,citecolor=blue,urlcolor=blue]{hyperref}

\doi{10.1214/154957804100000000}
\pubyear{0000}
\volume{0}
\firstpage{0}
\lastpage{0}

\startlocaldefs

\usepackage[utf8]{inputenc}
\usepackage{makeidx}
\usepackage{amsfonts}
\usepackage{amssymb}
\usepackage{gensymb}
\usepackage{caption}
\usepackage[labelformat=simple]{subcaption}
\usepackage{bbm}
\usepackage{lipsum}
\usepackage{enumitem}
\usepackage[english]{babel} 
\usepackage{bm} 
\usepackage[T1]{fontenc}	 
\usepackage{lmodern}
\usepackage{amsthm} 
\usepackage{graphicx} 
\usepackage{booktabs} 
\usepackage{multicol} 
\usepackage{tocstyle} 
\usepackage[toc,page]{appendix}

\newtheorem{theorem}{Theorem}[section]

\newtheorem{lemma}[theorem]{Lemma}
\usepackage[top=1in, bottom=1in, left=1.25in, right=1.25in]{geometry}
\DeclarePairedDelimiter\abs{\lvert}{\rvert}
\newcommand\bigfrac[2]{%
	\begin{array}{c}
		#1 \\
		\hline
		#2
	\end{array}}
	
\allowdisplaybreaks

\endlocaldefs
	
\doi{10.1214/154957804100000000}
\pubyear{0000}
\volume{0}
\firstpage{0}
\lastpage{0}

\begin{document}

\begin{frontmatter}

\title{Trace class Markov chains for the Normal-Gamma Bayesian shrinkage model}

\runtitle{Trace class Markov chains for Normal-Gamma model}

\author{\fnms{Liyuan} \snm{Zhang}\ead[label=e1]{zhangliyuan@stat.ufl.edu}}
\and
\author{\fnms{Kshitij} \snm{Khare}\ead[label=e2]{kdkhare@stat.ufl.edu}}
\address{}

\thankstext{t2}{Kshitij Khare (\printead{e2}) is Associate Professor, Department of 
Statistics, University of Florida. Liyuan Zhang (\printead{e1}) is Graduate 
Student, Department of Statistics, University of Florida.}

\affiliation{affl}{University of Florida}

\runauthor{Zhang and Khare}

\begin{abstract}
High-dimensional data, where the number of variables exceeds or is comparable to 
the sample size, is now pervasive in many scientific applications. In recent years, 
Bayesian shrinkage models have been developed as effective and computationally 
feasible tools to analyze such data, especially in the context of linear regression. 
In this paper, we focus on the Normal-Gamma shrinkage model developed by 
Griffin and Brown \cite{Griffin:Brown:2010}. This model subsumes the popular 
Bayesian lasso model, and a three-block Gibbs sampling algorithm to sample 
from the resulting intractable posterior distribution has been developed in 
\cite{Griffin:Brown:2010}. We consider an alternative two-block Gibbs sampling 
algorithm, and rigorously demonstrate its advantage over the three-block 
sampler by comparing specific spectral properties. In particular, we show that the 
Markov operator corresponding to the two-block sampler is trace class (and 
hence Hilbert-Schmidt), whereas the operator corresponding to the three-block 
sampler is not even Hilbert-Schmidt. The trace class property for the two-block 
sampler implies geometric convergence for the associated Markov chain, which 
justifies the use of Markov chain CLT's to obtain practical error bounds for MCMC 
based estimates. Additionally, it facilitates theoretical comparisons of the 
two-block sampler with sandwich algorithms which aim to improve performance 
by inserting inexpensive extra steps in between the two conditional draws of the 
two-block sampler. 
\end{abstract}




\tableofcontents

\end{frontmatter}

\section{Introduction}

\noindent
In recent years, the explosion of data, due to advances in science and information technology, 
has left almost no field untouched. The availability of high-throughput data from genomic, 
finance, environmental, marketing (among other) applications has created an urgent need for 
methodology and tools for analyzing high-dimensional data. In particular, consider the 
linear model $\bm{\mathrm{Y}} = X \bm{\beta} + \sigma{\boldsymbol \epsilon}$, where 
$\bm{\mathrm{Y}}$ is an $n \times 1$ real valued response variable, $X$ is a known $n 
\times p$ matrix, $\bm{\beta}$ is an unknown $p \times 1$ vector of regression coefficients, 
$\sigma$ is an unknown scale parameter and the entries of ${\boldsymbol \epsilon}$ are 
independent standard normals. In the high-dimensional datasets mentioned above, often 
$n < p$, and classical least squares methods fail. The lasso \cite{Tibshirani:1996} was 
developed to provide sparse estimates of the regression coefficient vector $\bm{\beta}$ in 
these sample-starved settings (several adaptations/alternatives have been proposed since 
then). It was observed in \cite{Tibshirani:1996} that the lasso estimate 
is the posterior mode obtained when one puts i.i.d Laplace priors on the elements of 
$\bm{\beta}$ (conditional on $\sigma$). This observation has led to a flurry of recent 
research concerning the development of prior distributions for $(\bm{\beta}, \sigma)$ 
that yield posterior distributions with high (posterior) probability around sparse values 
of $\bm{\beta}$, i.e., values of $\bm{\beta}$ that have many entries equal to $0$. Such 
prior distributions are referred to as ``continuous shrinkage priors " and the 
corresponding models are referred to as  ``Bayesian shrinkage models". 
Bayesian shrinkage methods have gained popularity and have been extensively used 
in a variety of applications including ecology, finance, image processing and neuroscience 
(see, for example, \cite{Yi:2008, deloscampos2009, Demiguel:2009, jacquemin2014, 
xing2012, mishchenko2012, gu2013, PongWong:2014, PongWong:Woolliams:2014}). 

In this paper, we focus on the well-known Normal-Gamma shrinkage model introduced 
in Griffin and Brown \cite{Griffin:Brown:2010}. The model is specified as follows: 
\begin{eqnarray}
& & \bm{\mathrm{Y}}~|~\bm{\beta},\bm{\tau},\sigma^{2} \sim {\mathcal{N}}_n \left( 
X\beta,\sigma^2I_n \right) \nonumber\\
& & \bm{\beta}~|~\sigma^2,\bm{\tau} \sim {\mathcal{N}}_p\left(\bm{0}_p,
\sigma^2D_{\bm{\tau}}\right)\nonumber\\
& & \sigma^2 \sim \mathrm{Inverse-Gamma}\left(\alpha,\xi\right)~ 
\left(\mathrm{allow~ for ~impropriety ~via} ~\alpha=0 ~\mathrm{or} ~\xi=0\right)
\nonumber\\
& & \tau_j \stackrel{i.i.d}{\sim} \mathrm{Gamma}(a,b)~~ \mathrm{for}~ 
j=1,2,...,p, \label{eq:model}
\end{eqnarray}

\noindent
where $\mathcal{N}_p$ denotes the $p$-variate normal density, and 
$D_{\bm{\tau}}$ is a diagonal matrix with diagonal entries given by 
$\{\tau_j\}^p_{j=1}$. Also, Inverse-Gamma$(\alpha, \xi)$ and Gamma$(a,b)$ denote 
the Inverse-Gamma and Gamma densities with shape parameters $\alpha$ and $a$, 
and rate parameters $\xi$ and $b$ respectively. The marginal density of $\bm{\beta}$ 
given $\sigma^2$ in the above model is given by 
$$
\pi(\bm{\beta} \mid \sigma^2) = \prod_{j=1}^p \frac{b^a}{\Gamma(a) \sqrt{2 \pi} \sigma} 
\left( \frac{\beta_j^2}{2b \sigma^2} \right)^{a/2} K_a \left( \frac{|\beta_j| \sqrt{2b}}
{\sigma} \right), 
$$

\noindent
where $K_a$ is the modified Bessel function of the second kind. The popular 
Bayesian lasso model of Park and Casella \cite{Park:Casella:2008} is a special 
case of the Normal-Gamma model above with $a = 1$, where the marginal 
density of $\bm{\beta}$ simplifies to 
$$
\pi(\bm{\beta} \mid \sigma^2) = \prod_{j=1}^p \frac{\sqrt{b}}{\sqrt{2} \sigma} 
e^{-\frac{|\beta_j| \sqrt{2b}}{\sigma}}. 
$$

\noindent
In this case, the marginal density for each $\beta_j$ (given $\sigma^2$) is the double 
exponential density. The Normal-Gamma family offers a wider choice for the tail 
behavior (as $a$ decreases, the marginal distribution becomes more peaked at zero, 
but has heavier tails), and thereby a more flexible mechanism for model shrinkage. 

The posterior density of $(\bm{\beta}, \sigma^2)$ for the Normal-Gamma model is 
intractable in the sense that closed form computation or direct sampling is not feasible. 
Griffin and Brown \cite{Griffin:Brown:2010} note that the full conditional densities 
of $\bm{\beta}, \sigma^2$ and $\bm{\tau}^2$ are easy to sample from, and develop 
a three-block Gibbs sampling Markov chain to generate samples from the 
desired posterior density. This Markov chain, denoted by $\tilde{\Phi} := 
\{(\tilde{\bm{\beta}}_m,\tilde{\sigma}_m^2)\}_{m=0}^{\infty}$ (on the state 
space $\mathbb{R}^p \times \mathbb{R}_{+}$), is driven by the Markov transition 
density (Mtd) 
\begin{equation} \label{eq:mtdthreeblock}
	\begin{split}
	\widetilde{k}\left(\left(\bm{\beta},\sigma^{2}\right),\left(\bar{\bm{\beta}},
	\bar{\sigma}^{2}\right)\right)
	&=\int_{\mathbb{R}_{+}^p}\pi\left(\bar{\sigma}^{2}\mid\bar{\bm{\beta}},\bm{\tau},
	\bm{\mathrm{Y}}\right)\pi\left(\bar{\bm{\beta}}\mid\sigma^{2},\bm{\tau},
	\bm{\mathrm{Y}}\right)\pi\left(\bm{\tau}\mid{\bm{\beta}},\sigma^{2},
	\bm{\mathrm{Y}}\right) d \bm{\tau}. 
	\end{split}
\end{equation}

\noindent
Here $\pi (\cdot \mid \cdot)$ denotes the conditional density of the first group of 
arguments given the second group of arguments. The one-step dynamics of this 
Markov chain to move from the current state, $\left( \tilde{\bm{\beta}}_m, 
\tilde{\sigma}^2_m \right)$, to the next state, $\left( \tilde{\bm{\beta}}_{m+1}, 
\tilde{\sigma}^2_{m+1} \right)$ can be described as follows: 
\begin{itemize}
\item Draw $\bm{\tau}$ from $\pi(\cdot \mid \tilde{\bm{\beta}}_m, \tilde{\sigma}^2_m, 
\bm{\mathrm{Y}})$. 
\item Draw $\tilde{\bm{\beta}}_{m+1}$ from $\pi(\cdot \mid \tilde{\sigma}^2_m, 
\bm{\tau}, \bm{\mathrm{Y}})$. 
\item Draw $\tilde{\sigma}^2_{m+1} $ from $\pi(\cdot \mid \tilde{\bm{\beta}}_{m+1}, 
\bm{\tau}, \bm{\mathrm{Y}})$. 
\end{itemize}

\noindent
In \cite{Pal:Khare:2014}, the authors show that the distribution of the Markov chain 
$\tilde{\Phi}$ converges to the desired posterior distribution at a geometric rate (as the 
number of steps converges to $\infty$). 

As mentioned previously, the Bayesian lasso Markov chain of \cite{Park:Casella:2008} 
is a special case of the Normal-Gamma Markov chain when $a = 1$. In recent work 
\cite{SRKZ:2016}, it was shown that a two-block version of the Bayesian lasso chain 
(and a variety of other chains for Bayesian regression) can be developed. The authors 
in \cite{SRKZ:2016} then focus their theoretical investigations on the Bayesian lasso, 
and show that the two-block Bayesian lasso chain has a better behaved spectrum 
than the original three-block Bayesian lasso chain in the following sense: the Markov 
operator corresponding to the two-block Bayesian lasso chain is trace class 
(eigenvalues are countable and summable, and hence in particular 
square-summable), while the Markov operator corresponding to the original three-block 
Bayesian lasso chain is not Hilbert-Schmidt (the corresponding absolute value 
operator either does not have a countable spectrum, or has a countable set of 
eigenvalues that are not square-summable). 

Based on the method outlined in \cite{SRKZ:2016}, a two-block version of the 
three-block Normal-Gamma Markov chain $\tilde{\Phi}$ can be constructed as 
follows. The two-block Markov chain, denoted by $\Phi = \{(\bm{\beta}_m,
\sigma_m^2)\}_{m=0}^{\infty}$ (on the state space $\mathbb{R}^p \times 
\mathbb{R}_{+}$), is driven by the Markov transition density (Mtd)
\begin{equation} \label{eqtwoblock}
k\left((\bm{\beta}, \sigma^2), (\tilde{\bm{\beta}}, \tilde{\sigma}^2)\right) = 
\int_{\mathbb{R}_{+}^p} \pi(\bm{\tau} \mid \bm{\beta}, \sigma^2, {\bf Y}) ~
\pi(\tilde{\bm{\beta}},\tilde{\sigma}^2 \mid \bm{\tau}, {\bf Y}) d \bm{\tau}. 
\end{equation}

\noindent
The one-step dynamics of this Markov chain to move from the current state, 
$\left( \bm{\beta}_m, \sigma^2_m \right)$, to the next state, $\left( \bm{\beta}_{m+1}, 
\sigma^2_{m+1} \right)$ can be described as follows: 
\begin{itemize}
\item Draw $\bm{\tau}$ from $\pi(\cdot \mid \bm{\beta}_m, \sigma^2_m, 
\bm{\mathrm{Y}})$. 
\item Draw $\sigma^2_{m+1} $ from $\pi(\cdot \mid \bm{\tau}, \bm{\mathrm{Y}})$ 
and draw $\bm{\beta}_{m+1}$ from $\pi(\cdot \mid \sigma^2_{m+1}, \bm{\tau}, 
\bm{\mathrm{Y}})$. 
\end{itemize}

\noindent
The goal of this paper is to investigate whether the theoretical results for the Bayesian 
lasso in \cite{SRKZ:2016} hold for the more general and complex setting of 
the Normal-Gamma model. In particular, we establish that the Markov operator 
corresponding to the two-block chain $\Phi$ is trace class when $a > \frac{1}{2}$ 
(Theorem \ref{thm1}). On the other hand, the Markov operator corresponding to the 
three-block chain $\tilde{\Phi}$ is not Hilbert-Schmidt for all values of $a$ 
(Theorem \ref{thm2}). These results hold {\bf for all} values of the sample size $n$ and 
the number of independent variables $p$. Since the Bayesian lasso is a special case 
with $a = 1$, our results subsume the spectral results in \cite{SRKZ:2016}. We note 
that establishing the results in the Normal-Gamma setup is {\bf much harder} than the 
Bayesian lasso setting. This is in part due to the fact that the modified Bessel function 
$K_a$ does not in general have a closed form when $a \neq 1$, and a heavy dose of 
various identities and bounds for these Bessel functions is needed to analyze the 
appropriate integrals in the proofs of Theorem \ref{thm1} and Theorem \ref{thm2}. 

We now discuss further some of the implications of establishing the trace-class 
property for $\Phi$. Note that the Markov chain $\Phi$ arises from a two-block Data 
Augmentation (DA) algorithm, with $(\bm{\beta}, \sigma^2)$ as the parameter 
block of interest and $\bm{\tau}$ as the augmented parameter block. Hence the 
corresponding Markov operator, denoted by $K$, is a positive, self-adjoint operator 
(see \cite{Hobert:Marchev:2008}). Establishing that a positive self-adjoint operator is 
trace class implies that it has a discrete spectrum, and that (countably many, 
non-negative) eigenvalues are summable. The trace class property implies 
compactness, which further implies geometric ergodicity of the underlying Markov 
chain (see \cite[Section 2]{PKH:2016}, for example). Geometric ergodicity, in turn, 
facilitates use of Markov chain central limit theorems to provide error bounds for 
Markov chain based estimates of relevant posterior expectations. The DA 
interpretation of $\Phi$ also enables us to use the Haar PX-DA technique 
from \cite{Hobert:Marchev:2008} and construct a ``sandwich" Markov chain 
by adding an inexpensive extra step in between the two conditional draws 
involved in one step of $\Phi$ (see Section \ref{sec:Haarsandwichchain} for details). 
The trace class property for $\Phi$, along with results in \cite{Khare:Hobert:2011}, 
implies that the sandwich chain is also trace class, and that each ordered eigenvalue 
of the sandwich chain is dominated by the corresponding ordered eigenvalue of 
$\Phi$ (with at least one strict domination). Recent work in 
\cite{Adamczak:Bednorz:2015} provides a rigorous approach to approximate 
eigenvalues of trace class Markov chains whose Mtd is available in closed form. 
These results are not applicable to the two-block sampler as its is not available in 
closed form, and extending results in \cite{Adamczak:Bednorz:2015} to such settings 
is a topic of ongoing research. 

The rest of the paper is organized as follows. In Section \ref{sec:form:densities}, we
provide the form of the relevant conditional densities for the Markov chains $\Phi$ 
and $\tilde{\Phi}$. In Section \ref{sec:twoblockDA}, we establish the trace class 
property for the two-block Markov chain $\Phi$. In Section \ref{sec:threeblockchain}, 
we show that the three-block Markov chain is not Hilbert-Schmidt. In 
Section \ref{sec:Haarsandwichchain}, we derive the Haar PX-DA sandwich chain 
corresponding to the two-block DA chain. Finally, in Section \ref{sec:experiments} 
we compare the performance of the two-block, three-block and the Haar PX-DA 
based chains on simulated and real datasets.

\section{Form of relevant densities} \label{sec:form:densities}

\noindent
In this section, we present expressions for various densities corresponding 
to the Normal-Gamma model in (\ref{eq:model}). These densities appear in the 
Mtd for the Markov chains $\Phi$ and $\tilde{\Phi}$. 

The joint density for the parameter vector $(\bm{\beta}, \bm{\tau} , \sigma^2)$ 
conditioned on the data vector $\bf y$ is given by the following:
\begin{equation} \label{eq:joint}
	\pi(\bm{\beta}, \bm{\tau} , \sigma^2| {\bf Y}) \propto
	\frac {e^{-\frac{\left(\left(\bm{\mathrm{Y}}-X\bm{\beta}\right)^T\left(\bm{\mathrm{Y}}-X\bm{\beta}\right)\right)}{2\sigma^2}}}{(\sqrt{2\pi})^n\sigma^n
	}\frac {e^{-\frac{\bm{\beta}^TD_{\bm{\tau}}^{-1}\bm{\beta}}{2\sigma^2}}}{(\sqrt{2\pi})^p\sigma^p}
	\left(\prod_{j=1}^{p}\tau_j^{a-\frac 1 2 -1}e^{-b\tau_j}\right)
	(\sigma^2)^{-a-1}e^{-\frac{\xi}{\sigma^2}}
\end{equation}
	
\noindent
Based on the joint density in (\ref{eq:joint}), the following conditional
distributions can be derived in a straightforward fashion.
	\vskip  0.2cm
\begin{itemize}
\item  $~~\bm{\beta}~|~\sigma^2,\bm{\tau},\bm{\mathrm{Y}}\sim \bm{\mathcal{N}_p}\left(\left(X^TX+D_{\bm{\tau}}^{-1}\right)^{-1}X^T\bm{\mathrm y},\sigma^2\left(X^TX+D_{\bm{\tau}}^{-1}\right)^{-1}\right)$
\end{itemize}
	
	In particular,
	\begin{equation} \label{eq:cond1}
	\pi\left(\bm{\beta}~|~\bm{\tau},\sigma^2,\bm{\mathrm{Y}}\right) = \frac{~|~X^TX+D^{-1}_{\bm{\tau}}~|~^{\frac{1}{2}}}{\sqrt{2\pi\sigma^2}^p}e^{-\frac{\left(\bm{\beta}-\left(X^TX+D^{-1}_{\bm{\tau}}\right)^{-1}X^T\bm{\mathrm{Y}}\right)^T\left(X^TX+D^{-1}_{\bm{\tau}}\right)\left(\bm{\beta}-\left(X^TX+D^{-1}_{\bm{\tau}}\right)^{-1}X^T\bm{\mathrm{Y}}\right)}{2\sigma^2}},
	\end{equation}
	for $\bm{\beta} \in \mathbb{R}^p.$

\begin{itemize}
\item  $~~\sigma^2~|~\bm{\beta},\bm{\tau},\bm{\mathrm{Y}}\sim \mathrm{Inverse-Gamma}\left(\frac{n+p+2\alpha}{2},\frac{\left(\bm{\mathrm{y}-X\bm{\beta}}\right)^T\left(\bm{\mathrm{y}-X\bm{\beta}}\right)+\bm{\beta}^TD_{\bm{\tau}}^{-1}\bm{\beta}+2\xi}{2}\right).$
\end{itemize}
	
In particular,
		\begin{equation} \label{eq:cond2}
		\begin{split}
		\pi\left(\sigma^2~|~\bm{\beta},\bm{\tau},\bm{\mathrm{Y}}\right) &= \frac{\left(\left(\bm{\mathrm{Y}}-X\bm{\beta}\right)^T\left(\bm{\mathrm{Y}}-X\bm{\beta}\right)+\bm{\beta}^TD_{\bm{\tau}}^{-1}\bm{\beta}+2\xi\right)^{\frac{n+p+2\alpha}{2}}}{2^{\frac{n+p+2\alpha}{2}}\Gamma\left(\frac{n+p+2\alpha}{2}\right)}\left(\sigma^2\right)^{-\frac{n+p+2\alpha}{2}-1}\\
		&{\times}~ e^{-\frac{\left(\left(\bm{\mathrm{Y}}-X\bm{\beta}\right)^T\left(\bm{\mathrm{Y}}-X\bm{\beta}\right)+\bm{\beta}^TD_{\bm{\tau}}^{-1}\bm{\beta}+2\xi\right)}{2\sigma^2}},\\
		\end{split}
		\end{equation}
		for $\sigma^2 \in \mathbb{R}_{+}.$
		
\begin{itemize}
\item  $~~\sigma^2~|~\bm{\tau},\bm{\mathrm{Y}}\sim \mathrm{Inverse-Gamma}\left(\frac{n+2\alpha}{2}, \frac{{{\bf Y}}^T\left(I-X{\bf A}_{\bm{\tau}}^{-1}X^T\right){{\bf Y}}+2\xi}{2}\right).$
\end{itemize}
	
In particular,
		\begin{equation} \label{eq:cond3}
		\begin{split}
		\pi\left(\sigma^2~|~\bm{\tau},\bm{\mathrm{Y}}\right) &= 
		\frac{	\left({{\bf Y}}^T\left(I-X{\bf A}_{\bm{\tau}}^{-1}X^T\right){{\bf Y}}+2\xi\right)^{\frac{n+2\alpha}{2}}}{2^{\frac{n+2\alpha}{2}}\Gamma\left(\frac{n+2\alpha}{2}\right)}
			\left(\sigma^2\right)^{-\frac{n+2\alpha}{2}-1}\\
			 & \times \exp \left( -\frac{1}{2\sigma^2}\left({{\bf Y}}^T\left(I-X{\bf A}_{\bm{\tau}}^{-1}X^T\right){{\bf Y}}+2\xi\right) \right) 
		,\\
		\end{split}
		\end{equation}
		for $\sigma^2 \in \mathbb{R}_{+}.$
		
\begin{itemize}
\item Given $\bm{\beta},\sigma^2$ and $\bf y,$ the variables $\tau_1,\tau_2,...,\tau_p$ 
are conditionally independent, and the conditional density of $\tau_j$ given $
\bm{\beta},\sigma^2$ and $\bf y$ is GIG$(a-\frac 1 2, 2b, \frac{\beta_j^2}{\sigma^2} ).$ 
\end{itemize}

In particular, 
	\begin{equation} \label{eq:cond4}
	\begin{split}
	\pi(\bm{\tau}\mid \bm{\beta},\sigma^2,{\bf Y}) = \prod_{j=1}^{p} 
	\frac{(2b\sigma^2)^{\frac {a-\frac 1 2}{2}}}{2\abs*{\beta_j}^{a-\frac 1 2}K_{a-\frac 1 2}\left(\sqrt{2b\frac{\beta_j^2}{\sigma^2}}\right) } \tau_j^{(a-\frac 1 2)-1}e^{-\frac 1 2\left\{2b\tau_j+\frac{\beta_j^2}{\sigma^2}\frac{1}{\tau_j}\right\}},
	\end{split}
	\end{equation}
for $\bm{\tau} \in \mathbb{R}_{+}^p.$

\section{Properties of the two-block Gibbs sampler} \label{sec:twoblockDA}

\noindent
In this section, we show that the operator associated with the two-block Gibbs sampler 
$\Phi,$ with Markov transition density $k$ specified in (\ref{eqtwoblock}) is trace class 
when $a>\frac{1}{2}$ and is not trace class when $0 < a \leq \frac{1}{2}$. 
\newtheorem{thm}{Theorem}
\begin{thm} \label{thm1}
For all values of $n$ and $p$, the Markov operator corresponding to the two-block 
Markov chain $\Phi$ is trace class (and hence Hilbert-Schmidt) when $a>\frac{1}{2}$ 
and is not trace class when $0 < a \leq \frac{1}{2}.$ 
\end{thm}

\noindent 
{\it Proof} In the current setting, the trace class property is equivalent to the finiteness 
of the integral (see \cite[Section 2]{PKH:2016}, for example)
\begin{equation} \label{eqtrace1}
\iint\limits_{\mathbb{R}^{p} \times \mathbb{R}_{+}} k \left(\left(\bm{\beta},\sigma^{2}\right), 
\left(\bm{\beta},\sigma^{2}\right)\right)\,d{\bm{\beta}}\,d{\sigma^2}. 
\end{equation}

\noindent
We will consider five separate cases: $a > 1$, $3/4 \leq a \leq 1$, $1/2 < a < 3/4$, 
$0 < a < 1/2$ and $a = 1/2$. In the first three cases, we will show that the integral 
in (\ref{eqtrace1}) is finite, and in the last two cases we will show that the integral 
in (\ref{eqtrace1}) is infinite. The proof is a lengthy and intricate algebraic exercise 
involving careful upper/lower bounds for modified Bessel functions and conditional 
densities, and we will try to provide a road-map/explanation whenever possible. 
We will start with the case $a > 1$. 

\medskip

\noindent
{\bf Case 1: $a > 1$} 

\noindent
By the definition of $k$, we have 
\begin{eqnarray}
& & \iint\limits_{\mathbb{R}_{+}^{p} \times \mathbb{R}_{+}}  k \left(\left(\bm{\beta},
\sigma^{2}\right), \left(\bm{\beta},\sigma^{2}\right)\right)\,d{\bm{\beta}}\,d{\sigma^2} 
\nonumber\\
&=& \iiint\limits_{\mathbb{R}_{+}^p \times \mathbb{R}_{+}^{p}\times \mathbb{R}_{+}}
\pi\left(\bm{\tau}~|~{\bm{\beta}},\sigma^{2},{\bf Y}\right)\pi\left(\bm{\beta},
\sigma^2\mid\bm{\tau},{\bf Y}\right)\,d{\bm{\beta}}\,d{\bm{\tau}}\,d{\sigma^2} 
\nonumber\\
&=& \iiint\limits_{\mathbb{R}_{+}^p \times \mathbb{R}_{+}^{p}\times \mathbb{R}_{+}}
\pi\left(\bm{\tau}~|~{\bm{\beta}},\sigma^{2},{\bf Y}\right)\pi\left(\sigma^{2}~|~\bm{\tau},
{\bf Y}\right)\pi\left(\bm{\beta}~|~\sigma^{2},\bm{\tau},{\bf Y}\right)\,d{\bm{\beta}}
\,d{\bm{\tau}}\,d{\sigma^2}. \label{eq:cond5}
\end{eqnarray}

\noindent
As a first step, we will gather all the terms with ${\boldsymbol \tau}$, and then focus 
on finding an upper bound for the inner integral with respect to ${\boldsymbol \tau}$. 
Using (\ref{eq:cond1}), (\ref{eq:cond3}) and (\ref{eq:cond4}), we get, 
\begin{eqnarray*}
& & \iint\limits_{\mathbb{R}_{+}^{p} \times \mathbb{R}_{+}}  k \left(\left(\bm{\beta},
\sigma^{2}\right), \left(\bm{\beta},\sigma^{2}\right)\right)\,d{\bm{\beta}}\,d{\sigma^2}\\
&=& C_1\iiint\limits_{\mathbb{R}^p \times \mathbb{R}_{+}^{p}\times \mathbb{R}_{+}}
\prod\limits_{j=1}^p\frac{\left(2b\sigma^2\right)^{\frac{a-\frac 1 2}{2}}}
{2\abs*{\beta_j}^{a-\frac 1 2}K_{a-\frac 1 2}\left(\sqrt{2b\frac{\beta_j^2}{\sigma^2}}
\right)}\tau_j^{\left(a-\frac 1 2\right)-1}e^{-\frac 1 2 \left\{2b\tau_j+\frac{\beta_j^2}
{\sigma^2}\frac{1}{\tau_j}\right\}} \left(\sigma^2\right)^{-\frac{n+2\alpha}{2}-1} 
\end{eqnarray*}

\begin{eqnarray}
& &  \times \exp \left( -\frac{1}{2\sigma^2}\left({{\bf Y}}^T\left(I-X{\bf A}_{\bm{\tau}}^{-1}X^T\right){{\bf Y}}+2\xi\right) \right) 
	\left({{\bf Y}}^T\left(I-X{\bf A}_{\bm{\tau}}^{-1}X^T\right){{\bf Y}}+2\xi\right)^{\frac{n+2\alpha}{2}}(\sigma^2)^{-\frac{p}{2}} \nonumber\\ 
& & \times |{\bf A}_{\bm{\tau}}~|~^{\frac{1}{2}} {\times} \exp \left( -\frac{\left(\bm{\beta}-{\bf A}_{\bm{\tau}}^{-1}X^T 
		{\bf Y}\right)^T{\bf A}_{\bm{\tau}}\left(\bm{\beta}-{\bf A}_{\bm{\tau}}^{-1}X^T{\bf Y}\right)}{2\sigma^2} \right) \,d{\bm{\beta}}
	\,d{\bm{\tau}}\,d{\sigma^2} \nonumber\\
	&\stackrel{(a)}{\leq}& C_2 \iiint\limits_{\mathbb{R}^p \times \mathbb{R}_{+}^{p}\times \mathbb{R}_{+}}
	\frac{\exp \left(-\frac{\xi}{\sigma^2} \right)}{\left(\sigma^2\right)^{\frac{n+2\alpha}{2}+1}}\frac{~|~{\bf A}_{\bm{\tau}}~|~^\frac{1}{2}}{\left(\sigma^2\right)^\frac{p}{2}}\prod\limits_{j=1}^p\frac{\left(2b\sigma^2\right)^{\frac{a-\frac 1 2}{2}}}{2\abs*{\beta_j}^{a-\frac 1 2}K_{a-\frac 1 2}\left(\sqrt{2b\frac{\beta_j^2}{\sigma^2}}\right)}\tau_j^{\left(a-\frac 1 2\right)-1}e^{-\frac 1 2 \left\{2b\tau_j+\frac{\beta_j^2}{\sigma^2}\frac{1}{\tau_j}\right\}} 
\nonumber\\
& & \times~{\exp \left(-\frac{1}{2\sigma^2}\left(\bm{\beta}^T{\bf A}_{\bm{\tau}}\bm{\beta}-2\bm{\beta}^TX^T{\bf Y}+{{\bf Y}}^T{\bf Y}\right) \right)}
	\,d{\bm{\beta}}\,d{\bm{\tau}}\,d{\sigma^2} \nonumber\\
&\stackrel{(a')}{\leq}& C_2 \iiint\limits_{\mathbb{R}^p \times \mathbb{R}_{+}^{p}\times \mathbb{R}_{+}}
	\frac{\exp \left(-\frac{\xi}{\sigma^2} \right)}{\left(\sigma^2\right)^{\frac{n+p+2\alpha}{2}+1}} {~|~{\bf A}_{\bm{\tau}}~|~^\frac{1}{2}} \prod\limits_{j=1}^p\frac{\left(2b\sigma^2\right)^{\frac{a-\frac 1 2}{2}}}{2\abs*{\beta_j}^{a-\frac 1 2}K_{a-\frac 1 2}\left(\sqrt{2b\frac{\beta_j^2}{\sigma^2}}\right)}\tau_j^{\left(a-\frac 1 2\right)-1}e^{-\frac 1 2 \left\{2b\tau_j+\frac{2 \beta_j^2}{\sigma^2}\frac{1}{\tau_j}\right\} }	\,d{\bm{\beta}}\,d{\bm{\tau}}\,d{\sigma^2},\nonumber
	\\ \label{eq:cond6}
\end{eqnarray}

\noindent
where 
$C_1=\frac{1}{(2\pi)^{\frac p 2}2^{\frac {n+2\alpha} 2}\Gamma\left(\frac{n+2\alpha}
{2}\right)}$ and $C_2 =  {\left({{\bf Y}}^T{{\bf Y}}+2\xi\right)^{\frac{n+2\alpha}{2}}}C_1$.
Note that $(a)$ follows from 
$$
{{\bf Y}}^T\left(I-X{\bf A}_{\bm{\tau}}^{-1}X^T\right){\bf Y} +2\xi \leq {{\bf Y}^T}{{\bf Y}}
+2\xi, 
$$

\noindent
and $(a')$ follows from 
\begin{align*}
&~~~~{\exp \left(-\frac{1}{2\sigma^2}\left(\bm{\beta}^T{\bf A}_{\bm{\tau}}
\bm{\beta}-2\bm{\beta}^TX^T{\bf Y}+{{\bf Y}}^T{\bf Y}\right) \right)} \\
& = \exp\left(-\frac{\bm{\beta}D_{\bm{\tau}}^{-1}\bm{\beta}}{2\sigma^2}\right)
{\exp \left\{-\frac{1}{2\sigma^2}\left(\bm{\beta}^TX^TX\bm{\beta}-2\bm{\beta}^TX^T 
{\bf Y}+{{\bf Y}}^T{\bf Y}\right) \right\}}\\
& \leq \exp\left(-\frac{\bm{\beta}D_{\bm{\tau}}^{-1}\bm{\beta}}{2\sigma^2}\right). 
\end{align*}

\noindent
We now focus on the inner integral in (\ref{eq:cond6}) defined by 
\begin{equation} \label{eqtrace2}
H\left(\bm{\beta},\sigma^2\right) \stackrel{\Delta}{=} 
\int_{\mathbb{R}_{+}^p}{~|~{\bf A}_{\bm{\tau}}~|~^\frac{1}{2}} 
\prod_{j=1}^p\frac{\left(2b\sigma^2\right)^{\frac{a-\frac 1 2}{2}}}{2\abs*{\beta_j}
^{a-\frac 1 2}K_{a-\frac 1 2}\left(\sqrt{2b\frac{\beta_j^2}{\sigma^2}}\right)}\tau_j^{\left(a-
\frac 1 2\right)-1}e^{-\frac 1 2 \left\{2b\tau_j+\frac{2 \beta_j^2}{\sigma^2}\frac{1}{\tau_j}
\right\} }d\bm{\tau} 
\end{equation}

\noindent
Let $\lambda$ denote the largest eigenvalue of $X^T X$. Using the definition of 
${\bf A}_{\bm{\tau}}$, it follows that 
\begin{eqnarray}
& & {~|~{\bf A}_{\bm{\tau}}~|~^\frac{1}{2}}
\prod\limits_{j=1}^p\frac{\left(2b\sigma^2\right)^{\frac{a-\frac 1 2}{2}}}{2\abs*{\beta_j}^{a-\frac 1 2}K_{a-\frac 1 2}\left(\sqrt{2b\frac{\beta_j^2}{\sigma^2}}\right)}\tau_j^{\left(a-
\frac 1 2\right)-1}e^{-\frac 1 2 \left\{2b\tau_j+\frac{2 \beta_j^2}{\sigma^2}\frac{1}{\tau_j}
\right\} } \nonumber\\
&\leq& \prod\limits_{j=1}^p \left(\lambda+\frac{1}{\sqrt{\tau_j}}\right)
\prod\limits_{j=1}^p\frac{\left(2b\sigma^2\right)^{\frac{a-\frac 1 2}{2}}}{2\abs*{\beta_j}^{a-\frac 1 2}K_{a-\frac 1 2}\left(\sqrt{2b\frac{\beta_j^2}{\sigma^2}}\right)}\tau_j^{\left(a-\frac 1 2\right)-1}e^{-\frac 1 2 \left\{2b\tau_j+\frac{2 \beta_j^2}{\sigma^2}\frac{1}{\tau_j}\right\} } \nonumber\\
&=& \left[\lambda^p+\left(\frac{1}{\tau_1}+...+\frac{1}{\sqrt{\tau_p}}\right)\lambda^{p-1}+
\left(\frac{1}{\sqrt{\tau_1\tau_2}}+...+\frac{1}{\sqrt{\tau_i\tau_j}}+...\right)
\lambda^{p-2}+...+\frac{1}{\sqrt{\tau_1\tau_2...\tau_p}}\right]\prod\limits_{j=1}^p c_j \nonumber\\
\label{eq:cond7}
\end{eqnarray}

where 
$$
c_j = \frac{\left(2b\sigma^2\right)^{\frac{a-\frac 1 2}{2}}}{2\abs*{\beta_j}^{a-\frac 1 2}
K_{a-\frac 1 2}\left(\sqrt{2b\frac{\beta_j^2}{\sigma^2}}\right)}\tau_j^{\left(a-\frac 1 
2\right)-1}e^{-\frac 1 2 \left\{2b\tau_j+\frac{2\beta_j^2}{\sigma^2}\frac{1}{\tau_j}\right\}}. 
$$ 

\noindent
We now examine a generic term of the sum in (\ref{eq:cond7}). Note that $c_j$ and 
$\frac{c_j}{\sqrt{\tau_j}}$ are both (unnormalized) GIG densities. Hence, for any 
subset $\mathcal{L} = \{\ell_1, \ell_2, \cdots, \ell_m\}$ of $\{1,2, \cdots, p\}$, using 
the form of the GIG density, we get 
\begin{eqnarray}
& & \int_{\mathbb{R}_{+}^p} \frac{1}{\sqrt{\tau_{\ell_1} \tau_{\ell_2} \cdots 
\tau_{\ell_m}}} \prod\limits_{j=1}^p c_j d\bm{\tau} \nonumber\\
&=& \left( \prod\limits_{j \notin \mathcal{L}} \int c_j d{\tau_j} \right) \times \left( \prod\limits_{j \in \mathcal{L}} \int \frac{c_j}{\sqrt{\tau_j}} d{\tau_j} \right) 
\nonumber\\
&=& \left( \prod\limits_{j \notin \mathcal{L}} (\sqrt{2})^{a-\frac{1}{2}}\frac{K_{a-\frac 1 2}\left(\sqrt{4b\frac{\beta_j^2}{\sigma^2}}\right)}{K_{a-\frac 1 2}\left(\sqrt{2b\frac{\beta_j^2}{\sigma^2}}\right)} \right) \times \left( \prod\limits_{j \in \mathcal{L}} 
\left(2b\sigma^2\right)^{\frac 1 4}\left(\sqrt{2}\right)^{a-1}\abs*{\beta_{j}}^{-\frac 1 2}~\frac{K_{a- 1}\left(\sqrt{4b\frac{\beta_{j}^2}{\sigma^2}}\right)}{K_{a-\frac 1 2}\left(\sqrt{2b\frac{\beta_{j}^2}{\sigma^2}}\right)} \right). \label{eq:cond8}
\end{eqnarray}

\noindent First, by \cite[Page 266]{Laforgia:1991}, we get that 
$$
\frac{K_{a-\frac 1 2}\left(\sqrt{4b\frac{\beta_j^2}{\sigma^2}}\right)}{K_{a-\frac 1 
2}\left(\sqrt{2b\frac{\beta_j^2}{\sigma^2}}\right)} < \exp\left({\sqrt{2b\frac{\beta_j^2}
{\sigma^2}}-\sqrt{4b\frac{\beta_j^2}{\sigma^2}}}\right) < \exp \left( -\frac{\sqrt{b}\abs*{\beta_j}}
{2 \sigma} \right)
$$

\noindent
for all $a > \frac{1}{2}$. Next, using the fact that if $x>0,$ then $\nu \rightarrow K_{\nu}
(x)$ is an increasing function for $\nu>0$ (again, see \cite[Page 266]{Laforgia:1991}), 
we get 
$$
\frac{K_{a- 1}\left(\sqrt{4b\frac{\beta_j^2}{\sigma^2}}\right)}{K_{a-\frac 1 2}\left(\sqrt{2b
\frac{\beta_j^2}{\sigma^2}}\right)} \leq \frac{K_{a-\frac 1 2}\left(\sqrt{4b\frac{\beta_j^2}
{\sigma^2}}\right)}{K_{a-\frac 1 2}\left(\sqrt{2b\frac{\beta_j^2}{\sigma^2}}\right)} < 
\exp \left( -\frac{\sqrt{b}\abs*{\beta_j}}{2 \sigma} \right). 
$$

\noindent
Hence, from (\ref{eq:cond8}), we get that 
\begin{equation} \label{eq:cond9}
\int_{\mathbb{R}_{+}^p} \frac{1}{\sqrt{\tau_{\ell_1} \tau_{\ell_2} \cdots 
\tau_{\ell_m}}} \prod\limits_{j=1}^p c_j d\bm{\tau} < \sqrt{2}^{pa}b^{\frac{\abs*{\mathcal{L}}}{4}}(\sigma^2)^{\frac{\abs*{\mathcal{L}}}{4}} \left( \prod_{j \in 
\mathcal{L}} \abs*{\beta_j}^{-\frac 1 2} \right) \exp \left( -\sum_{j=1}^p 
\frac{\sqrt{b}\abs*{\beta_j}}{2 \sigma} \right). 
\end{equation}

\noindent
It follows from (\ref{eqtrace2}), (\ref{eq:cond7}) and (\ref{eq:cond9}) that 
$$
H\left(\bm{\beta},\sigma^2\right) \leq \sum_{\mathcal{L} \subseteq \{1,2, \cdots, p\}} 
\lambda^{p-|\mathcal{L}|} \sqrt{2}^{pa} b^{\frac{\abs*{\mathcal{L}}}{4}}(\sigma^2)^{\frac{\abs*{\mathcal{L}}}{4}}\left( \prod_{j \in 
\mathcal{L}} \abs*{\beta_j}^{-\frac 1 2} \right) \exp \left( -\sum_{j=1}^p 
\frac{\sqrt{b}\abs*{\beta_j}}{2 \sigma} \right). 
$$

\noindent
By (\ref{eq:cond6}) and (\ref{eqtrace2}), the trace class property will be established if 
we show that for every $\mathcal{L} \subseteq \{1,2, \cdots, p\}$, the integral 
$$
\iint\limits_{\mathbb{R}_{+} \times \mathbb{R}^p} \frac{\exp \left(-\frac{\xi}{\sigma^2} \right)}{\left(\sigma^2\right)^{\frac{n+p+2\alpha}{2}-\frac{\abs*{\mathcal{L}}}{4}+1}} \left( \prod_{j \in 
\mathcal{L}} \abs*{\beta_j}^{-\frac 1 2} \right) \exp \left( -\sum_{j=1}^p 
\frac{\sqrt{b}\abs*{\beta_j}}{2 \sigma} \right) \,d{\bm{\beta}}\,d{\sigma^2} 
$$

\noindent
is finite. We proceed to show this by first simplifying the inner integral with respect to 
$\bm{\beta}$. Using the form of the Gamma density, we get 
\begin{eqnarray}
& & \int_{\mathbb{R}^p} \left( \prod_{j \in \mathcal{L}} \abs*{\beta_j}^{-\frac 1 2} \right) 
\exp \left( -\sum_{j=1}^p \frac{\sqrt{b}\abs*{\beta_j}}{2 \sigma} \right) \,d{\bm{\beta}} 
\nonumber\\
&=& \left( \prod_{j \notin \mathcal{L}} \int_\mathbb{R} \exp \left( -\frac{\sqrt{b}\abs*{\beta_j}}{2 
\sigma} \right) \,d{\beta_j} \right) \times \left( \prod_{j \in \mathcal{L}} \int_\mathbb{R} 
\abs*{\beta_j}^{-\frac 1 2} \exp \left( -\frac{\sqrt{b}\abs*{\beta_j}}{2 \sigma} \right) \,d{\beta_j} 
\right) \nonumber\\
&=& \left(\frac{4 \sigma}{\sqrt{b}}\right)^{p-|\mathcal{L}|} \left( 2 \Gamma \left( \frac{1}{2} \right) \sqrt{\frac{2 \sigma}{\sqrt{b}}} 
\right)^{|\mathcal{L}|} \nonumber\\
&\leq& \frac{8^p}{(\sqrt{b})^{p-\frac{|\mathcal{L}|}{2}}} \sigma^{p - \frac{|\mathcal{L}|}{2}}. \label{eqtrace3}
\end{eqnarray}

\noindent
It follows by (\ref{eqtrace3}) that 
\begin{eqnarray*}
& & \iint\limits_{\mathbb{R}_{+} \times \mathbb{R}^p} \frac{\exp \left(-\frac{\xi}
{\sigma^2} \right)}{\left(\sigma^2\right)^{\frac{n+p+2\alpha}{2}-\frac{\abs*{\mathcal{L}}}{4} +1}} \left( \prod_{j \in 
\mathcal{L}} \abs*{\beta_j}^{-\frac 1 2} \right) \exp \left( -\sum_{j=1}^p 
\frac{\sqrt{b}\abs*{\beta_j}}{2 \sigma} \right) \,d{\bm{\beta}}\,d{\sigma^2}\\
&\leq& \frac{8^p}{(\sqrt{b})^{p-\frac{|\mathcal{L}|}{2}}}\int_{\mathbb{R}_+} \frac{\exp \left(-\frac{\xi}{\sigma^2} \right)}
{\left(\sigma^2\right)^{\frac{n+p+2\alpha}{2}-\frac{\abs*{\mathcal{L}}}{4}+1}} \sigma^{p - \frac{|\mathcal{L}|}{2}} 
\,d{\sigma^2}\\
&\leq& \frac{8^p}{(\sqrt{b})^{p-\frac{|\mathcal{L}|}{2}}} \int_{\mathbb{R}_+} \frac{\exp \left(-\frac{\xi}{\sigma^2} \right)}
{\left(\sigma^2\right)^{\frac{n}{2}+\alpha+1}} \,d{\sigma^2}\\
&=& \frac{8^p \Gamma \left( \frac{n}{2} + \alpha \right)}
{{(\sqrt{b})^{p-\frac{|\mathcal{L}|}{2}}}\xi^{\frac{n}{2} + \alpha}} < \infty. 
\end{eqnarray*}

\noindent
As discussed above, this establishes the trace class property in the case $a > 1$. 
 
\medskip

\noindent
{\bf Case 2: $3/4 \leq a \leq 1$}

\noindent
In this case, we first note that all arguments in Case 1 go through verbatim until 
(\ref{eq:cond8}). Next, we note that 
\begin{align}
\frac{K_{a- 1}\left(\sqrt{4b\frac{\beta_j^2}{\sigma^2}}\right)}{K_{a-\frac 1 2}\left(\sqrt{2b\frac{\beta_j^2}{\sigma^2}}\right)}=
{\frac{K_{a-\frac 1 2}\left(\sqrt{4b\frac{\beta_j^2}{\sigma^2}}\right)}{K_{a-\frac 1 2}\left(\sqrt{2b\frac{\beta_j^2}{\sigma^2}}\right)}}~
{\frac{K_{a- 1}\left(\sqrt{4b\frac{\beta_j^2}{\sigma^2}}\right)}{K_{a-\frac 1 2}\left(\sqrt{4b\frac{\beta_j^2}{\sigma^2}}\right)}}\label{case2a}
\end{align}

\noindent If $a\geq \frac 3 4,$ then $ a-\frac{1}{2}>0$, and by 
\cite[Page 266]{Laforgia:1991}, we get 
\begin{align}
\frac{K_{a-\frac 1 2}\left(\sqrt{4b\frac{\beta_j^2}{\sigma^2}}\right)}{K_{a-\frac 1 
		2}\left(\sqrt{2b\frac{\beta_j^2}{\sigma^2}}\right)} < \exp\left({\sqrt{2b\frac{\beta_j^2}
		{\sigma^2}}-\sqrt{4b\frac{\beta_j^2}{\sigma^2}}}\right) < \exp \left( -\frac{\sqrt{b}\abs*{\beta_j}}
{2 \sigma} \right) \label{case2b}
\end{align}
Using the property that $K_{\nu}(x)=K_{-\nu}(x)$ (see \cite{Abramowitz:Stegun:1965}, Page 375), we obtain
\begin{align*}
{\frac{K_{a- 1}\left(\sqrt{4b\frac{\beta_j^2}{\sigma^2}}\right)}{K_{a-\frac 1 2}\left(\sqrt{4b\frac{\beta_j^2}{\sigma^2}}\right)}} = {\frac{K_{1- a}\left(\sqrt{4b\frac{\beta_j^2}{\sigma^2}}\right)}{K_{a-\frac 1 2}\left(\sqrt{4b\frac{\beta_j^2}{\sigma^2}}\right)}}
\end{align*}
If $ \frac 3 4 \leq a <1,$  then $0 < 1-a \leq a-\frac 1 2$. Since $\nu \rightarrow K_{\nu}
(x)$ is increasing in $\nu>0$ for $x > 0$ (see \cite[Page 266]{Laforgia:1991}), it follows that ${K_{1- a}\left(\sqrt{4b\frac{\beta_j^2}{\sigma^2}}\right)}\leq {K_{a-\frac 1 2}\left(\sqrt{4b\frac{\beta_j^2}{\sigma^2}}\right)}$ for $\frac{3}{4} \leq a < 1$. 
Also, by the integral formula (see \cite{Abramowitz:Stegun:1965}, Page 376)
\[K_{\nu}(t)= \int_{0}^{\infty} e^{-t\cosh z} \cosh(\nu z) \,dz, \nu\in\mathcal{R}. 
\]
Since $\cosh(\nu z) \geq \cosh(0)$ for any $\nu>0, z>0$ ($x \rightarrow \cosh(x)$ is 
increasing on $[0,\infty)$), we get \[K_{\nu}(t) \geq K_0(t)\] for $\nu>0$. In particular, 
$K_0(\sqrt{4b\frac{\beta_j^2}{\sigma^2}}) \leq K_{\frac 1 2}(\sqrt{4b\frac{\beta_j^2}
{\sigma^2}}).$ Hence for all $a \in \left[\frac 3 4 , 1\right],$ we have
\begin{align}
{\frac{K_{a- 1}\left(\sqrt{4b\frac{\beta_j^2}{\sigma^2}}\right)}{K_{a-\frac 1 2}\left(\sqrt{4b
\frac{\beta_j^2}{\sigma^2}}\right)}} \leq 1. 
\label{case2c}
\end{align}

\noindent
It follows from (\ref{case2a}), (\ref{case2b}) and (\ref{case2c}) that 
\[\frac{K_{a- 1}\left(\sqrt{4b\frac{\beta_j^2}{\sigma^2}}\right)}{K_{a-\frac 1 2} 
\left(\sqrt{2b\frac{\beta_j^2}{\sigma^2}}\right)} < \exp \left( -\frac{\sqrt{b}\abs*{\beta_j}}
{2 \sigma} \right) \]

\noindent
Now, using exactly the same arguments as in the proof of Case 1 (following 
(\ref{eq:cond8})) the trace class property can be shown the case $\frac{3}{4} \leq a 
\leq 1$. 
 
\medskip

\noindent
{\bf Case 3: $1/2 < a < 3/4$}

\noindent
Again, in this case, we first note that all arguments in Case 1 go through verbatim until 
(\ref{eq:cond8}). Also, by \cite[Page 266]{Laforgia:1991} and $K_\nu (x) = K_{-\nu} (x)$ 
for $x > 0$, we get 
\begin{align*}
\frac{K_{a- 1}\left(\sqrt{4b\frac{\beta_j^2}{\sigma^2}}\right)}{K_{a-\frac 1 2}\left(\sqrt{2b\frac{\beta_j^2}{\sigma^2}}\right)}&=
{\frac{K_{a- 1}\left(\sqrt{4b\frac{\beta_j^2}{\sigma^2}}\right)}{K_{a-\frac 1 2}\left(\sqrt{4b\frac{\beta_j^2}{\sigma^2}}\right)}}
\frac{K_{a-\frac 1 2}\left(\sqrt{4b\frac{\beta_j^2}{\sigma^2}}\right)}{K_{a-\frac 1 2}\left(\sqrt{2b\frac{\beta_j^2}{\sigma^2}}\right)} < \exp \left( -\frac{\sqrt{b}\abs*{\beta_j}}
{2 \sigma} \right) 
{\frac{K_{1- a}\left(\sqrt{4b\frac{\beta_j^2}{\sigma^2}}\right)}{K_{a-\frac 1 2}\left(\sqrt{4b\frac{\beta_j^2}{\sigma^2}}\right)}}
\end{align*}

\noindent
Note that if $1/2 < a < 3/4$, then $1-a-(a-\frac 1 2 ) = \frac 3 2 -2a \in (0,\frac 1 2)$.  
It follows by \cite[Page 640]{Pal:Khare:2014}) that   
\begin{align*}
\frac{K_{1- a}\left(\sqrt{4b\frac{\beta_j^2}{\sigma^2}}\right)}{K_{a-\frac 1 2}\left(\sqrt{4b\frac{\beta_j^2}{\sigma^2}}\right)} \leq 
\dfrac{\left(2a\right)^{\frac 3 2 -2a}}{\left(\sqrt{4b}\frac{\abs*{\beta_j}}{\sigma}\right)^{\frac 3 2 - 2a}} +1 
\end{align*}
Hence,
\[\frac{K_{a- 1}\left(\sqrt{4b\frac{\beta_j^2}{\sigma^2}}\right)}{K_{a-\frac 1 2}\left(\sqrt{2b\frac{\beta_j^2}{\sigma^2}}\right)} < \exp \left( -\frac{\sqrt{b}\abs*{\beta_j}}
{2 \sigma} \right) \left( \dfrac{\left(2-2a\right)^{\frac 3 2 -2a}}{\left(\sqrt{4b}\frac{\abs*{\beta_j}}{\sigma}\right)^{\frac 3 2 - 2a}} +1 \right).\]

\noindent
By (\ref{eq:cond8}), for any subset $\mathcal{L} = \{\ell_1, \ell_2, \cdots, \ell_m\}$ of 
$\{1,2, \cdots, p\}$ we get 
\begin{eqnarray}
& & \int_{\mathbb{R}_{+}^p} \frac{1}{\sqrt{\tau_{\ell_1} \tau_{\ell_2} \cdots 
		\tau_{\ell_m}}} \prod\limits_{j=1}^p c_j d\bm{\tau} \nonumber\\
&=& \left( \prod\limits_{j \notin \mathcal{L}} (\sqrt{2})^{a-\frac{1}{2}}\frac{K_{a-\frac 1 2}\left(\sqrt{4b\frac{\beta_j^2}{\sigma^2}}\right)}{K_{a-\frac 1 2}\left(\sqrt{2b\frac{\beta_j^2}{\sigma^2}}\right)} \right) \times \left( \prod\limits_{j \in \mathcal{L}} 
\left(2b\sigma^2\right)^{\frac 1 4}\left(\sqrt{2}\right)^{a-1}\abs*{\beta_{j}}^{-\frac 1 2}~\frac{K_{a- 1}\left(\sqrt{4b\frac{\beta_{j}^2}{\sigma^2}}\right)}{K_{a-\frac 1 2}\left(\sqrt{2b\frac{\beta_{j}^2}{\sigma^2}}\right)} \right)\nonumber\\
&\leq&  \sqrt{2}^{pa}b^{\frac{\abs*{\mathcal{L}}}{4}}(\sigma^2)^{\frac{\abs*{\mathcal{L}}}{4}} \left( \prod_{j \in 
	\mathcal{L}} \abs*{\beta_j}^{-\frac 1 2} \right) \prod_{j\in \mathcal{L}}\left( \dfrac{\left(2-2a\right)^{\frac 3 2 -2a}}{\left(\sqrt{4b}\frac{\abs*{\beta_j}}{\sigma}\right)^{\frac 3 2 - 2a}} +1 \right) \exp \left( -\sum_{j=1}^p 
\frac{\sqrt{b}\abs*{\beta_j}}{2 \sigma} \right). \label{eqtrace4}
\end{eqnarray}

\noindent
It follows from (\ref{eqtrace2}) that 
$$
H\left(\bm{\beta},\sigma^2\right) \leq \sum_{\mathcal{L} \subseteq \{1,2, \cdots, p\}} 
\lambda^{p-|\mathcal{L}|}\sqrt{2}^{pa}b^{\frac{\abs*{\mathcal{L}}}{4}}(\sigma^2)^{\frac{\abs*{\mathcal{L}}}{4}} \left( \prod_{j \in 
	\mathcal{L}} \abs*{\beta_j}^{-\frac 1 2} \right) \prod_{j\in \mathcal{L}}\left( \dfrac{\left(2-2a\right)^{\frac 3 2 -2a}}{\left(\sqrt{4b}\frac{\abs*{\beta_j}}{\sigma}\right)^{\frac 3 2 - 2a}} +1 \right) \exp \left( -\sum_{j=1}^p 
\frac{\sqrt{b}\abs*{\beta_j}}{2 \sigma} \right)$$

\noindent
By (\ref{eq:cond6}), the trace class property will be established if we show that for 
every ${\mathcal{L} \subseteq \{1,2, \cdots, p\}} , $ the integral 
\begin{align}
\iint\limits_{\mathbb{R}_{+} \times \mathbb{R}^p} \frac{\exp \left(-\frac{\xi}{\sigma^2} 
\right)}{\left(\sigma^2\right)^{\frac{n+p+2\alpha}{2}-\frac{\abs*{\mathcal{L}}}{4}+1}} 
\left( \prod_{j \in \mathcal{L}} \abs*{\beta_j}^{-\frac 1 2} \right) \prod_{j\in \mathcal{L}}
\left( \dfrac{\left(2-2a\right)^{\frac 3 2 -2a}}{\left(\sqrt{4b}\frac{\abs*{\beta_j}}{\sigma}
\right)^{\frac 3 2 - 2a}} +1 \right) \exp \left( -\sum_{j=1}^p \frac{\sqrt{b}\abs*{\beta_j}}{2 
\sigma} \right) \,d{\bm{\beta}}\,d{\sigma^2} \label{eqtrace5}
\end{align}

\noindent
is finite. We proceed to show this by first integrating out $\bm{\beta}$. Using the form 
of the Gamma density, we get 
\begin{align}
&~~~~~\int\limits_{\mathbb{R}^p} \left( \prod_{j \in 
	\mathcal{L}} \abs*{\beta_j}^{-\frac 1 2} \right) \prod_{j\in \mathcal{L}}\left( \dfrac{\left(2-2a\right)^{\frac 3 2 -2a}}{\left(\sqrt{4b}\frac{\abs*{\beta_j}}{\sigma}\right)^{\frac 3 2 - 2a}} +1 \right) \exp \left( -\sum_{j=1}^p 
\frac{\sqrt{b}\abs*{\beta_j}}{2 \sigma} \right) \,d{\bm{\beta}}\nonumber\\
&=\left( \prod_{j \notin \mathcal{L}} \int_\mathbb{R} \exp \left( -\frac{\sqrt{b}\abs*{\beta_j}}{2 
	\sigma} \right) \,d{\beta_j} \right) \times \left( \prod_{j \in \mathcal{L}} \int_\mathbb{R} ~
\abs*{\beta_j}^{-\frac 1 2}\left( \dfrac{\left(2-2a\right)^{\frac 3 2 -2a}}{\left(\sqrt{4b}\frac{\abs*{\beta_j}}{\sigma}\right)^{\frac 3 2 - 2a}} +1 \right) \exp \left( -\frac{\sqrt{b}\abs*{\beta_j}}{2 \sigma} \right) \,d{\beta_j} 
\right) \nonumber\\
&=\left(\frac{4 \sigma}{\sqrt{b}}\right)^{p-|\mathcal{L}|}\prod_{j \in \mathcal{L}} \left(\int\limits_{\mathbb{R}}\abs*{\beta_j}^{-\frac 1 2}\exp \left( -\frac{\sqrt{b}\abs*{\beta_j}}{2 \sigma} \right)\,d{\beta_j}     
+\int\limits_{\mathbb{R}}\abs*{\beta_j}^{-\frac 1 2}
\dfrac{\left(2-2a\right)^{\frac 3 2 -2a}}{\left(\sqrt{4b}\frac{\abs*{\beta_j}}{\sigma}\right)^{\frac 3 2 - 2a}}
\exp \left( -\frac{\sqrt{b}\abs*{\beta_j}}{2 \sigma} \right) 
\right)\,d{\beta_j}\nonumber\\
&= \left(\frac{4 \sigma}{\sqrt{b}}\right)^{p-|\mathcal{L}|}
\prod_{j \in \mathcal{L}} \left(2 \Gamma \left( \frac{1}{2} \right) \sqrt{\frac{2 \sigma}{\sqrt{b}}}+\left(1-a\right)^{\frac 3 2 -2a}\Gamma(2a-1)2^{2a}\left(\sqrt{b}\right)^{-\frac 1 2}\sigma^{\frac 1 2}
\right)\nonumber\\
&=\left(\frac{4 \sigma}{\sqrt{b}}\right)^{p-|\mathcal{L}|} 
\sigma^{\frac{|\mathcal{L}|}{2}}\left(\sqrt{b}\right)^{-\frac{|\mathcal{L}|}{2}}C_3(\mathcal{L})\nonumber\\
&\leq4^p\left(\sqrt{b}\right)^{\frac{|\mathcal{L}|}{2}-p}C_3(\mathcal{L})\sigma^{p-\frac{|\mathcal{L}|}{2}}\label{eqtrace6}
\end{align}

\noindent
where $C_3(\mathcal{L})=\left(2\sqrt{2} \Gamma(\frac{1}{2})+(1-a)^{\frac 3 2 -2a}
\Gamma(2a-1)2^{2a}\right)^{|\mathcal{L}|}$. It follows by (\ref{eqtrace5}) that 
 \begin{align*}
 &~~\iint\limits_{\mathbb{R}_{+} \times \mathbb{R}^p} \frac{\exp \left(-\frac{\xi}{\sigma^2} \right)}{\left(\sigma^2\right)^{\frac{n+p+2\alpha}{2}-\frac{\abs*{\mathcal{L}}}{4}+1}} \left( \prod_{j \in 
 	\mathcal{L}} \abs*{\beta_j}^{-\frac 1 2} \right) \prod_{j\in \mathcal{L}}\left( \dfrac{\left(2-2a\right)^{\frac 3 2 -2a}}{\left(\sqrt{4b}\frac{\abs*{\beta_j}}{\sigma}\right)^{\frac 3 2 - 2a}} +1 \right) \exp \left( -\sum_{j=1}^p 
 \frac{\sqrt{b}\abs*{\beta_j}}{2 \sigma} \right) \,d{\bm{\beta}}\,d{\sigma^2} \\
 &  \leq 4^p\left(\sqrt{b}\right)^{\frac{|\mathcal{L}|}{2}-p}C_3(\mathcal{L})\int_{\mathbb{R}_+} \frac{\exp \left(-\frac{\xi}{\sigma^2} \right)}
 {\left(\sigma^2\right)^{\frac{n}{2}+\alpha+1}} \,d{\sigma^2}\\
 &=\frac{4^p\left(\sqrt{b}\right)^{\frac{|\mathcal{L}|}{2}-p}C_3(\mathcal{L})\Gamma(\frac{n}{2}+\alpha)}{\xi^{\frac{n}{2}+\alpha}} < \infty.
 \end{align*}

\noindent
As discussed above, this establishes the trace class property in the case $1/2 < a < 
3/4$. 

\medskip

\noindent
{\bf Case 4: $0 <a < 1/2$}

\noindent
Now, we'll show that when $a \in \left( 0,\frac 1 2 \right)$, 

$$
\iint\limits_{\mathbb{R}_{+}^{p} \times \mathbb{R}_{+}}  k \left(\left(\bm{\beta},
\sigma^{2}\right), \left(\bm{\beta},\sigma^{2}\right)\right)\,d{\bm{\beta}}\,d{\sigma^2} = 
\infty. 
$$
 
\noindent
Note that 
\begin{align*}
\frac{K_{a- 1}\left(\sqrt{4b\frac{\beta_j^2}{\sigma^2}}\right)}{K_{a-\frac 1 2}\left(\sqrt{2b\frac{\beta_j^2}{\sigma^2}}\right)}
&=
{\frac{K_{1-a }\left(\sqrt{4b\frac{\beta_j^2}{\sigma^2}}\right)}{K_{\frac 1 2-a}\left(\sqrt{2b\frac{\beta_j^2}{\sigma^2}}\right)}}. 
\end{align*}

\noindent
By \cite[Page 375]{Abramowitz:Stegun:1965}, if $\nu>0,$ then $ \frac{K_{\nu}(x)}
{\frac{2^\nu\Gamma(\nu)}{2x^v}}\rightarrow 1$ as $x\rightarrow 0.$ Let 
$y={\sqrt{\frac{2b\beta_j^2}{\sigma^2}}}$. It follows that 
$\dfrac{\frac{K_{1-a }(\sqrt{2}y)}{\frac{2^{1-a}\Gamma(1-a)}{2(\sqrt{2}y)^{1-a}}}}
{\frac{{K_{\frac 1 2 -a}(y)}}{\frac{2^{\frac 1 2 -a}\Gamma(\frac 1 2 - a)}{~2y^{\frac 1 2 -
a}}}} \rightarrow 1$ as $y \rightarrow 0$. Hence there exists $\epsilon_1 > 0 $ such 
that $\dfrac{\frac{K_{1-a }(\sqrt{2}y)}{\frac{2^{1-a}\Gamma(1-a)}{2(\sqrt{2}y)^{1-a}}}}
{\frac{{K_{\frac 1 2 -a}(y)}}{\frac{2^{\frac 1 2 -a}\Gamma(\frac 1 2 - a)}{~2y^{\frac 1 2 -
a}}}} \geq \frac 1 2$ for $0 <y < \epsilon_1. $ Thus if $0 <y < \epsilon_1,$ we have 
\[{\frac{K_{1-a }\left(\sqrt{4b\frac{\beta_j^2}{\sigma^2}}\right)}{K_{\frac 1 2-a}
\left(\sqrt{2b\frac{\beta_j^2}{\sigma^2}}\right)}}= \frac{K_{1-a}(\sqrt{2}y)}{K_{\frac 1 2-a}
(y)} \geq \frac{\left(\sqrt{2}\right)^{a}\Gamma(1-a)}{2\Gamma(\frac 1 2 -a)} \frac{1}
{\sqrt{y}}
\]

\noindent
Since $K_{\nu}(x) >0 $ for positive $\nu$ and $x,$ we have 
\begin{eqnarray*}
\frac{K_{1-a}(\sqrt{2}y)}{K_{\frac 1 2-a}(\sqrt{2}y)} \geq \frac{K_{1-a}(\sqrt{2}y)}
{K_{\frac 1 2-a}(\sqrt{2}y)} I_{(0 <y < \epsilon_1)} \geq \frac{\left(\sqrt{2}\right)^{a}
\Gamma(1-a)}{2\Gamma(\frac 1 2 -a)} \frac{I_{(0 <y < \epsilon_1)}}{\sqrt{y}}, 
\end{eqnarray*}

\noindent
Using $y={\sqrt{\frac{2b\beta_j^2}{\sigma^2}}}$, we get 
\begin{align}
\frac{K_{a- 1}\left(\sqrt{4b\frac{\beta_j^2}{\sigma^2}}\right)}{K_{a-\frac 1 2}\left(\sqrt{2b
\frac{\beta_j^2}{\sigma^2}}\right)}\geq & \frac{\left(\sqrt{2}\right)^{a}\Gamma(1-a)}
{2\left(2b\right)^{\frac{1}{4}}\Gamma(\frac 1 2 -a)}
\frac{\sigma^\frac{1}{2}}{\abs*{\beta_j}^{\frac 1 2}}I_{(0 < \abs*{\beta_j}<\frac{\sigma
\epsilon_1}{\sqrt{2b}})} \label{eqtrace7}
\end{align}

\noindent
It follows from (\ref{eq:cond1}), (\ref{eq:cond3}) and (\ref{eq:cond4}) that 
\begin{align}
&  ~~~~~~~\iint\limits_{\mathbb{R}_{+}^{p} \times \mathbb{R}_{+}}  k \left(\left(\bm{\beta},\sigma^{2}\right), \left(\bm{\beta},\sigma^{2}\right)\right)\,d{\bm{\beta}}\,d{\sigma^2}\notag \\
&= C_1\iiint\limits_{\mathbb{R}^p \times \mathbb{R}_{+}^{p}\times \mathbb{R}_{+}}
\prod\limits_{j=1}^p\frac{\left(2b\sigma^2\right)^{\frac{a-\frac 1 2}{2}}}{2\abs*{\beta_j}^{a-\frac 1 2}K_{a-\frac 1 2}\left(\sqrt{2b\frac{\beta_j^2}{\sigma^2}}\right)}\tau_j^{\left(a-\frac 1 2\right)-1}e^{-\frac 1 2 \left\{2b\tau_j+\frac{\beta_j^2}{\sigma^2}\frac{1}{\tau_j}\right\} }
\left(\sigma^2\right)^{-\frac{n+2\alpha}{2}-1}\notag\\
&  \times \exp \left( -\frac{1}{2\sigma^2}\left({{\bf Y}}^T\left(I-X{\bf A}_{\bm{\tau}}^{-1}X^T\right){{\bf Y}}+2\xi\right) \right) 
\left({{\bf Y}}^T\left(I-X{\bf A}_{\bm{\tau}}^{-1}X^T\right){{\bf Y}}+2\xi\right)^{\frac{n+2\alpha}{2}}(\sigma^2)^{-\frac{p}{2}}\notag\\ 
&  \times |{\bf A}_{\bm{\tau}}~|~^{\frac{1}{2}} {\times} \exp \left( -\frac{\left(\bm{\beta}-{\bf A}_{\bm{\tau}}^{-1}X^T 
	{\bf Y}\right)^T{\bf A}_{\bm{\tau}}\left(\bm{\beta}-{\bf A}_{\bm{\tau}}^{-1}X^T{\bf Y}\right)}{2\sigma^2} \right) \,d{\bm{\beta}}
\,d{\bm{\tau}}\,d{\sigma^2}\label{eqtrace8}
\end{align}

\noindent
Furthermore, we have 
\begin{equation}
|{\bf A}_{\bm{\tau}}|^{\frac 1 2}= |X^TX+ D_{\bm{\tau}}^{-1}|^{\frac 1 2} \geq |
D_{\bm{\tau}}^{-1}|^{\frac 1 2},~~ ~~ \left({{\bf Y}}^T\left(I-X{\bf A}_{\bm{\tau}}^{-1}X^T
\right){{\bf Y}}\right)+2\xi \geq 2\xi, \label{eqtrace8.1}
\end{equation}

\noindent
and
\begin{align}
& ~~~~~\exp \left( -\frac{1}{2\sigma^2}\left({{\bf Y}}^T\left(I-X{\bf A}_{\bm{\tau}}^{-1}X^T
\right){{\bf Y}}+2\xi\right) \right) \exp \left( -\frac{\left(\bm{\beta}-{\bf A}_{\bm{\tau}}^{-1}
X^T {\bf Y}\right)^T{\bf A}_{\bm{\tau}}\left(\bm{\beta}-{\bf A}_{\bm{\tau}}^{-1}X^T{\bf Y}
\right)}{2\sigma^2} \right) \notag\\
& = e^{-\frac{2\xi}{2\sigma^2}} \exp \left\{
-\frac{1}{2\sigma^2} \left({{\bf Y}}^T\left(I-X{\bf A}_{\bm{\tau}}^{-1}X^T\right){{\bf Y}}+
{\left(\bm{\beta}-{\bf A}_{\bm{\tau}}^{-1}X^T {\bf Y}\right)^T{\bf A}_{\bm{\tau}}
\left(\bm{\beta}-{\bf A}_{\bm{\tau}}^{-1}X^T{\bf Y}\right)} \right) \right\}\notag\\
& = e^{-\frac{2\xi}{2\sigma^2}} \exp \left\{ -\frac{1}{2\sigma^2} \left({{\bf Y}}^T{{\bf Y}}
+ \bm{\beta}^T({X^TX+D^{-1}_{\bm{\tau}}})\bm{\beta}-2\bm{\beta}^TX^T{\bf Y} \right)
\right\} \notag\\
& = e^{-\frac{2\xi+{{\bf Y}}^T{{\bf Y}}}{2\sigma^2}} \exp\left(-\frac{\bm{\beta}
D_{\bm{\tau}}^{-1}\bm{\beta}}{2\sigma^2}\right) \exp \left\{ -\frac{1}{2\sigma^2} \left(
\bm{\beta}^TX^TX\bm{\beta}-2\bm{\beta}^TX^T{\bf Y} \right) \right\} \label{eqtrace8.2}
\end{align}

\noindent
If we denote the entries of $X^TX$ and $X^TY$ by $a_{ij}, b_i$ separately. It's easy to 
see there is at least $i$ such that $a_{ii}>0$ (if not, $a_{ii}=0 $ for all $i$, indicating 
$X$ is exactly $\bm{0}$.) Without loss of generality, we assume $a_{11}>0,$ so 
\begin{equation}
\exp \left\{ -\frac{1}{2\sigma^2} \left( \bm{\beta}^TX^TX\bm{\beta}-2\bm{\beta}^TX^T 
{\bf Y} \right) \right\} = g(\sigma^2,\beta_2,\cdots,\beta_p) \exp \left\{ -\frac{a_{11} 
\left( \beta_1+c\right)^2}{2\sigma^2} \right\} \label{eqtrace9}
\end{equation}

\noindent
where $g(\sigma^2,\beta_2,\cdots,\beta_p) = \exp\left\{\dfrac{a_{11}c^2}
{2\sigma^2}+2\sum\limits_{i=2}^{p}\beta_ib_i-\sum\limits_{2\leq i, j\leq p}a_{ij}\beta_i
\beta_j\right\}$ and $~c=\dfrac{a_{12}\beta_2+a_{13}\beta_3+\cdots+a_{1p}\beta_p-
b_1}{a_{11}}$. It follows from (\ref{eqtrace8}), (\ref{eqtrace8.1}), (\ref{eqtrace8.2}) 
and (\ref{eqtrace9}) that 
\begin{align}
&  ~~~~~~~~~\iint\limits_{\mathbb{R}_{+}^{p} \times \mathbb{R}_{+}}  k \left(\left(\bm{\beta},\sigma^{2}\right), \left(\bm{\beta},\sigma^{2}\right)\right)\,d{\bm{\beta}}\,d{\sigma^2}\notag \\
& \geq 
C_1\iiint\limits_{\mathbb{R}^p \times \mathbb{R}_{+}^{p}\times \mathbb{R}_{+}}
\frac{  e^{-\frac{2\xi+{{\bf Y}}^T{{\bf Y}}}{2\sigma^2}} }{\left(\sigma^2\right)^{\frac{n+2\alpha}{2}+1}}~(\sigma^2)^{-\frac{p}{2}}~
g(\sigma^2,\beta_2,\cdots,\beta_p) \exp \left\{
-\frac{a_{11}\left(
	\beta_1+c\right)^2}{2\sigma^2}
\right\}~|D_{\bm{\tau}}^{-1}|^{\frac 1 2}~ \exp\left(-\frac{\bm{\beta}D_{\bm{\tau}}^{-1}\bm{\beta}}{2\sigma^2}\right)~\notag \\
&
~~~\times \prod\limits_{j=1}^p\frac{\left(2b\sigma^2\right)^{\frac{a-\frac 1 2}{2}}}{2\abs*{\beta_j}^{a-\frac 1 2}K_{a-\frac 1 2}\left(\sqrt{2b\frac{\beta_j^2}{\sigma^2}}\right)}\tau_j^{\left(a-\frac 1 2\right)-1}e^{-\frac 1 2 \left\{2b\tau_j+\frac{\beta_j^2}{\sigma^2}\frac{1}{\tau_j}\right\} }
 \,d{\bm{\beta}}
\,d{\bm{\tau}}\,d{\sigma^2} \notag \\
& = 
{C_1}\iint\limits_{\mathbb{R}^p \times \mathbb{R}_{+}}
\frac{  e^{-\frac{2\xi+{{\bf Y}}^T{{\bf Y}}}{2\sigma^2}} }{\left(\sigma^2\right)^{\frac{n+p+2\alpha}{2}+1}}~~g(\sigma^2,{\beta_2},\cdots,\beta_p)~\exp \left\{
-\frac{a_{11}}{2\sigma^2}\left(
\beta_1+c\right)^2
\right\}
\int\limits_{\mathbb{R}_{+}^{p}}|D_{\bm{\tau}}^{-1}|^{\frac 1 2}~ \exp\left(-\frac{\bm{\beta}D_{\bm{\tau}}^{-1}\bm{\beta}}{2\sigma^2}\right)~\notag \\
&
~~~\times \prod\limits_{j=1}^p\frac{\left(2b\sigma^2\right)^{\frac{a-\frac 1 2}{2}}}{2\abs*{\beta_j}^{a-\frac 1 2}K_{a-\frac 1 2}\left(\sqrt{2b\frac{\beta_j^2}{\sigma^2}}\right)}\tau_j^{\left(a-\frac 1 2\right)-1}e^{-\frac 1 2 \left\{2b\tau_j+\frac{\beta_j^2}{\sigma^2}\frac{1}{\tau_j}\right\} }
\,d{\bm{\tau}}
\,d{\bm{\beta}}
\,d{\sigma^2}\label{greatcond}
\end{align}

\noindent
By (\ref{eqtrace7}), the inner integral can be bounded below as 
\begin{align}
&~~~~~\int\limits_{\mathbb{R}_{+}^{p}}|D_{\bm{\tau}}^{-1}|^{\frac 1 2}~ \exp\left(-\frac{\bm{\beta}D_{\bm{\tau}}^{-1}\bm{\beta}}{2\sigma^2}\right)~ \prod\limits_{j=1}^p\frac{\left(2b\sigma^2\right)^{\frac{a-\frac 1 2}{2}}}{2\abs*{\beta_j}^{a-\frac 1 2}K_{a-\frac 1 2}\left(\sqrt{2b\frac{\beta_j^2}{\sigma^2}}\right)}\tau_j^{\left(a-\frac 1 2\right)-1}e^{-\frac 1 2 \left\{2b\tau_j+\frac{\beta_j^2}{\sigma^2}\frac{1}{\tau_j}\right\} }\,d{\bm{\tau}}\notag \\
& =\prod\limits_{j=1}^p~~ \int\limits_{\mathbb{R}_{+}}
\frac{\left(2b\sigma^2\right)^{\frac{a-\frac 1 2}{2}}}{2\abs*{\beta_j}^{a-\frac 1 2}K_{a-\frac 1 2}\left(\sqrt{2b\frac{\beta_j^2}{\sigma^2}}\right)}\tau_j^{\left(a- 1 \right)-1}e^{-\frac 1 2 \left\{2b\tau_j+\frac{2\beta_j^2}{\sigma^2}\frac{1}{\tau_j}\right\} }\,d{{\tau_j}}\notag \\
& = \prod\limits_{j=1}^p
\left(2b\sigma^2\right)^{\frac 1 4}\left(\sqrt{2}\right)^{a-1}\abs*{\beta_{j}}^{-\frac 1 2}~\frac{K_{a- 1}\left(\sqrt{4b\frac{\beta_{j}^2}{\sigma^2}}\right)}{K_{a-\frac 1 2}\left(\sqrt{2b\frac{\beta_{j}^2}{\sigma^2}}\right)}\nonumber\\
& \geq
\prod\limits_{j=1}^p
\left(2b\sigma^2\right)^{\frac 1 4}\left(\sqrt{2}\right)^{a-1}\abs*{\beta_{j}}^{-\frac 1 2}~
\frac{\left(\sqrt{2}\right)^{a}\Gamma(1-a)}{2\left(2b\right)^{\frac{1}{4}}\Gamma(\frac 1 2 -a)}
\frac{\sigma^\frac{1}{2}}{\abs*{\beta_j}^{\frac 1 2}}I_{(0 < \abs*{\beta_j}<\frac{\sigma\epsilon}{\sqrt{2b}})} \notag
 \\ 
& = \frac{\left(\sqrt{2}\right)^{2a-3}\Gamma(1-a)}{\Gamma(\frac{1}{2}-a)} (\sigma^2)^{\frac p 2 } \prod\limits_{j=1}^p \frac{1}{{\abs*{\beta_j}}}I_{(0 < \abs*{\beta_j}<\frac{\sigma\epsilon}{\sqrt{2b}})}\label{innercond}
\end{align}

\noindent
It follows from (\ref{greatcond}) and (\ref{innercond}) that
\begin{align}
&  ~~~~\iint\limits_{\mathbb{R}_{+}^{p} \times \mathbb{R}_{+}}  k \left(\left(\bm{\beta},\sigma^{2}\right), \left(\bm{\beta},\sigma^{2}\right)\right)\,d{\bm{\beta}}\,d{\sigma^2} \notag\\
& \geq \frac{\left(\sqrt{2}\right)^{2a-3}\Gamma(1-a)C_1}{\Gamma(\frac{1}{2}-a)} \iint\limits_{\mathbb{R}^p \times \mathbb{R}_{+}}
\frac{  e^{-\frac{2\xi+{{\bf Y}}^T{{\bf Y}}}{2\sigma^2}} }{\left(\sigma^2\right)^{\frac{n+2\alpha}{2}+1}}~~g(\sigma^2,{\beta_2},\cdots,\beta_p)\exp \left\{
-\frac{a_{11}\left(
	\beta_1+c\right)^2}{2\sigma^2}
\right\} \prod\limits_{j=1}^p \frac{1}{{\abs*{\beta_j}}}I_{(0 < \abs*{\beta_j}<\frac{\sigma\epsilon}{\sqrt{2b}})}\,d{\bm{\beta}}
\,d{\sigma^2} \notag\\
& = C_1^* \int\limits_{\mathbb{R}_+} \frac{  e^{-\frac{2\xi+{{\bf Y}}^T{{\bf Y}}}
{2\sigma^2}} }{\left(\sigma^2\right)^{\frac{n+2\alpha}{2}+1}} \left\{ 
\int_{-\frac{\sigma\epsilon}{\sqrt{2b}}}^{\frac{\sigma\epsilon}{\sqrt{2b}}} \frac{1}
{\abs*{\beta_1}} \exp \left\{ -\frac{a_{11}\left( \beta_1+c\right)^2}{2\sigma^2}
\right\}\, d \beta_1 \right\} \times \notag\\
& ~~~~ \left\{ \int_{\mathbb{R}^{p-1}} g(\sigma^2,{\beta_2},\cdots,
\beta_p) \prod\limits_{j=2}^p \frac{1}{{\abs*{\beta_j}}}I_{(0 < \abs*{\beta_j}<\frac{\sigma
\epsilon}{\sqrt{2b}})}\,\prod\limits{j=2}^p d \beta_j \right\}\, d \sigma^2 \label{eqtrace10}
\end{align}

\noindent
where $C_1^* = \frac{\left(\sqrt{2}\right)^{2a-3}\Gamma(1-a)C_1}{\Gamma(\frac{1}{2}-
a)}$. However, we note that 
$$
\int_{-\frac{\sigma\epsilon}{\sqrt{2b}}}^{\frac{\sigma\epsilon}{\sqrt{2b}}} \frac{1}
{\abs*{\beta_1}} \exp \left\{ -\frac{a_{11}\left( \beta_1+c\right)^2}{2\sigma^2}
\right\}\, d \beta_1 \geq \int\limits_{0}^{\frac{\sigma\epsilon}{\sqrt{2b}}} \frac{1}
{\beta_1}\exp \left\{-\frac{a_{11}\left( \beta_1+c\right)^2}{2\sigma^2} \right\}\, 
d{\beta_1}=\infty 
$$

\noindent
where the last step follows from Propositon A1. By (\ref{eqtrace10}), it follows that the 
operator corresponding to the Markov transition density $k$ is not trace class when 
$0 < a < 1/2$. 

\medskip

\noindent
{\bf Case 5: $a= 1/2$}

\noindent
Finally, we show that when $a = \frac 1 2$, we have 
$$
\iint\limits_{\mathbb{R}_{+}^{p} \times \mathbb{R}_{+}}  k \left(\left(\bm{\beta},
\sigma^{2}\right), \left(\bm{\beta},\sigma^{2}\right)\right)\,d{\bm{\beta}}\,d{\sigma^2} = 
\infty. 
$$

\noindent
When $a =\frac 1 2, ~	\frac{K_{a- 1}\left(\sqrt{4b\frac{\beta_j^2}{\sigma^2}}\right)}{K_{a-\frac 1 2}\left(\sqrt{2b\frac{\beta_j^2}{\sigma^2}}\right)}=	\frac{K_{- \frac 1 2}\left(\sqrt{4b\frac{\beta_j^2}{\sigma^2}}\right)}{K_0\left(\sqrt{2b\frac{\beta_j^2}{\sigma^2}}\right)}=\frac{K_{\frac 1 2}\left(\sqrt{4b\frac{\beta_j^2}{\sigma^2}}\right)}{K_0\left(\sqrt{2b\frac{\beta_j^2}{\sigma^2}}\right)}.$ 
By \cite[Page 375]{Abramowitz:Stegun:1965}, if $z\rightarrow 0,$ then $K_0(z) \sim -\ln(z)$ and $K_{\frac 1 2}(z) \sim \frac{\Gamma(\frac 1 2)}{2}\left(\frac{2}{z}\right)^{\frac 1 2}$. As we did in Case 4, let $y={\sqrt{\frac{2b\beta_j^2}{\sigma^2}}}$. It follows that 
$$
\dfrac{\frac{K_{\frac 1 2}(\sqrt{2}y)}{\frac{\Gamma(\frac 1 2)}{2}\left(\frac{2}{\sqrt{2}y}\right)^{\frac 1 2}}}{\frac{K_0(y)}{-\ln(y)}}\rightarrow 1  ~~~~\text{as~} y \rightarrow 0. 
$$

\noindent
Hence there exists $\epsilon_2 \in (0,1)$ such that $\frac{K_{\frac 1 2}(\sqrt{2}y)}
{K_0(y)} \geq \frac{\Gamma(\frac{1}{2})2^{\frac 1 4}}{4}\frac{1}{-\sqrt{y}\ln(\sqrt{2}y)}$ 
for $0 < y < \epsilon_2$. It follows that 
\begin{align*}
{\frac{K_{\frac 1 2  }\left(\sqrt{4b\frac{\beta_j^2}{\sigma^2}}\right)}{K_0\left(\sqrt{2b
\frac{\beta_j^2}{\sigma^2}}\right)}}=~\frac{{K_{\frac 1 2}(\sqrt{2}y)}}{{K_0(y)}}&\geq 
\frac{{K_{\frac 1 2}(\sqrt{2}y)}}{{K_0(y)}}I_{\left(0<y<\epsilon_2\right)}\\ &\geq 
\frac{\Gamma(\frac{1}{2})2^{\frac 1 4}}{4}\frac{1}{-\sqrt{y}\ln(\sqrt{2}y)}I_{\left(0<y<
\epsilon_2\right)}\\ &=\frac{C_5(\sigma^2)^{\frac 1 4}}{{\abs*{\beta_j}^{\frac 1 2 }\left(-
\ln(\sqrt{2b})+\ln\sigma-\ln\abs*{\beta_j}\right)}}I_{(0 < \abs*{\beta_j}<\frac{\sigma
\epsilon_2}{\sqrt{2b}})}, 
\end{align*}
 
\noindent
where $C_5= \frac{\Gamma(\frac{1}{2})}{4b^{\frac 1 4}}$. We use this to get a 
lower bound for the inner integral with respect to $\bm{\tau}$ in (\ref{greatcond}). 
In particular, we note that 
\begin{align*}
	&~~~~\int\limits_{\mathbb{R}_{+}^{p}}|D_{\bm{\tau}}^{-1}|^{\frac 1 2}~ \exp\left(-\frac{\bm{\beta}D_{\bm{\tau}}^{-1}\bm{\beta}}{2\sigma^2}\right)~ \prod\limits_{j=1}^p\frac{\left(2b\sigma^2\right)^{\frac{a-\frac 1 2}{2}}}{2\abs*{\beta_j}^{a-\frac 1 2}K_{a-\frac 1 2}\left(\sqrt{2b\frac{\beta_j^2}{\sigma^2}}\right)}\tau_j^{\left(a-\frac 1 2\right)-1}e^{-\frac 1 2 \left\{2b\tau_j+\frac{\beta_j^2}{\sigma^2}\frac{1}{\tau_j}\right\} }\,d{\bm{\tau}}\notag \\
	& = \prod\limits_{j=1}^p
	\left(2b\sigma^2\right)^{\frac 1 4}\left(\sqrt{2}\right)^{a-1}\abs*{\beta_{j}}^{-\frac 1 2}~\frac{K_{-\frac 1 2 }\left(\sqrt{4b\frac{\beta_{j}^2}{\sigma^2}}\right)}{K_{0}\left(\sqrt{2b\frac{\beta_{j}^2}{\sigma^2}}\right)}\notag   \\
	& \geq
	\prod\limits_{j=1}^p
	\left(2b\sigma^2\right)^{\frac 1 4}\left(\sqrt{2}\right)^{a-1}\abs*{\beta_{j}}^{-\frac 1 2}~
	\frac{C_5(\sigma^2)^{\frac 1 4 }}{\abs*{\beta_j}^{\frac 1 2 }\left(-\ln(\sqrt{2b})+\ln\sigma-\ln\abs*{\beta_j}\right)}I_{(0 < \abs*{\beta_j}<\frac{\sigma\epsilon_2}{\sqrt{2b}})} \notag \\ 
	& =(2b)^{\frac p 4} (\sqrt{2})^{p(a-1)} C_5^p~(\sigma^2)^{\frac p 2} 	~\prod\limits_{j=1}^p \frac{I_{(0 < \abs*{\beta_j}<\frac{\sigma\epsilon_2}{\sqrt{2b}})}}{{\abs*{\beta_j}\left(-\ln(\sqrt{2b})+\ln\sigma-\ln\abs*{\beta_j}\right)}}
\end{align*}

\noindent
Using (\ref{greatcond}), it follows that 
\begin{align*}
	&  ~~~~~~~~~\iint\limits_{\mathbb{R}_{+}^{p} \times \mathbb{R}_{+}}  k \left(\left(\bm{\beta},\sigma^{2}\right), \left(\bm{\beta},\sigma^{2}\right)\right)\,d{\bm{\beta}}\,d{\sigma^2}\notag \\
	& \geq ~ 
{C_1}\iint\limits_{\mathbb{R}^p \times \mathbb{R}_{+}}
\frac{  e^{-\frac{2\xi+{{\bf Y}}^T{{\bf Y}}}{2\sigma^2}} }{\left(\sigma^2\right)^{\frac{n+p+2\alpha}{2}+1}}~~g(\sigma^2,{\beta_2},\cdots,\beta_p)~\exp \left\{
-\frac{a_{11}}{2\sigma^2}\left(
\beta_1+c\right)^2
\right\}
\int\limits_{\mathbb{R}_{+}^{p}}|D_{\bm{\tau}}^{-1}|^{\frac 1 2}~ \exp\left(-\frac{\bm{\beta}D_{\bm{\tau}}^{-1}\bm{\beta}}{2\sigma^2}\right) \times\\
& ~~~~ \times \prod\limits_{j=1}^p\frac{\left(2b\sigma^2\right)^{\frac{a-\frac 1 2}{2}}}{2\abs*{\beta_j}^{a-\frac 1 2}K_{a-\frac 1 2}\left(\sqrt{2b\frac{\beta_j^2}{\sigma^2}}\right)}\tau_j^{\left(a-\frac 1 2\right)-1}e^{-\frac 1 2 \left\{2b\tau_j+\frac{\beta_j^2}{\sigma^2}\frac{1}{\tau_j}\right\} }
\,d{\bm{\tau}}
\,d{\bm{\beta}}
\,d{\sigma^2}\\
& \geq (2b)^{\frac p 4} (\sqrt{2})^{p(a-1)} C_1C_5^p\iint\limits_{\mathbb{R}^{p-1} \times \mathbb{R}_{+}}
\frac{  e^{-\frac{2\xi+{{\bf Y}}^T{{\bf Y}}}{2\sigma^2}} }{\left(\sigma^2\right)^{\frac{n+2\alpha}{2}+1}} g(\sigma^2,{\beta_2},\cdots,\beta_p)~\prod\limits_{j=2}^p \frac{I_{(0 < \abs*{\beta_j}<\frac{\sigma\epsilon_2}{\sqrt{2b}})}}{{\abs*{\beta_j}\left(-\ln(\sqrt{2b})+\ln\sigma-\ln\abs*{\beta_j}\right)}} \times\\
& ~~~~ \left\{ \int\limits_{R} \exp \left\{ -\frac{a_{11}}{2\sigma^2}\left( \beta_1+c\right)^2
\right\} \frac{I_{(0 < \abs*{\beta_1}<\frac{\sigma\epsilon_2}{\sqrt{2b}})}}
{{\abs*{\beta_1}\left(-\ln(\sqrt{2b})+\ln\sigma-\ln\abs*{\beta_1}\right)}} \, d \beta_1 
\right\} \, d \bm{\beta}' \, d \sigma^2, 
\end{align*}

\noindent
where $\bm{\beta}'= (\beta_2,\cdots,\beta_p)$. By Proposition A2, we obtain 
\begin{align*}
	&~~~~~\int\limits_{R} \exp \left\{
	-\frac{a_{11}}{2\sigma^2}\left(
	\beta_1+c\right)^2
	\right\}
	\frac{I_{(0 < \abs*{\beta_1}<\frac{\sigma\epsilon_2}{\sqrt{2b}})}}{{\abs*{\beta_1}\left(-\ln(\sqrt{2b})+\ln\sigma-\ln\abs*{\beta_1}\right)}}
	\,d{\beta_1}\notag \\ & \geq 
	\int\limits_{0}^{\frac{\sigma\epsilon_2}{\sqrt{2b}}} 
	\frac{\exp \left\{
		-\frac{a_{11}}{2\sigma^2}\left(
		\beta_1+c\right)^2
		\right\}}{{{\beta_1}\left(-\ln(\sqrt{2b})+\ln\sigma-\ln{\beta_1}\right)}}
	\,d{\beta_1}\\
	 & = \infty 
\end{align*}

\noindent
It follows that the operator corresponding to the Markov transition density $k$ is not 
trace class when $a = \frac{1}{2}$. \hfill$\Box$

\section{Properties of the three-block Gibbs sampler} \label{sec:threeblockchain}

\noindent
In this section, we show that when $a>0, $ the Markov operator corresponding to the three-block Gibbs sampler $\tilde{\Phi},$ with Markov transition density $\tilde{k}$ specified in (1), is not Hilbert-Schmidt. Let $\tilde{K}$ be the Markov operator corresponding to $\tilde{\Phi}.$ We prove the following result. 

\begin{thm} \label{thm2}
	For all $a>0,$ the Markov operator $\tilde{K}$ is not Hilbert-Schmidt 
	{for all} possible values of $p$ and $n$. 
\end{thm}

\noindent
{\it Proof}~ Note that the Markov operator $\tilde{K}$ corresponding to the density 
$\tilde{k}$ is Hilbert-Schmidt if and only if $\tilde{K}^* \tilde{K}$ is trace class (see 
\cite{Jorgens:1982}, for example). Here $\tilde{K}^*$ denotes the adjoint of $\tilde{K}$. 
It follows that $\hat{K}$ is Hilbert-Schmidt if and only if $I < \infty$, where 
\begin{eqnarray}
I 
&:=& \int_{\mathbb{R}^p}\int_{\mathbb{R}_{+}} \int_{\mathbb{R}^p} \int_{\mathbb{R}_{+}}
\tilde{k} \bigg(\left(\bm{\beta},\sigma^{2}\right) , \left(\tilde{\bm{\beta}},\tilde{\sigma}^2\right)\bigg) 
\tilde{k}^* \bigg(\left(\bm{\beta},\sigma^{2}\right) , \left(\tilde{\bm{\beta}},\tilde{\sigma}^2\right)\bigg) d\bm{\beta}d\sigma^2
d\tilde{\bm{\beta}}d\tilde{\sigma^2} \nonumber\\
&=& \int_{\mathbb{R}^p}\int_{\mathbb{R}_{+}} \int_{\mathbb{R}^p} \int_{\mathbb{R}_{+}}
\tilde{k}^2\bigg(\left(\bm{\beta},\sigma^{2}\right)~,~\left(\tilde{\bm{\beta}},\tilde{\sigma}^2\right)\bigg) 
\frac{ \pi\left(\bm{\beta},{\sigma^{2}}~|~{\bf Y}\right)  }{ \pi\left(\tilde{\bm{\beta}},\tilde{\sigma}^2~|~{\bf Y}\right)  }d\bm{\beta} 
d\sigma^2d\tilde{\bm{\beta}}d\tilde{\sigma^2}. \label{eqtrace11}
\end{eqnarray}

\noindent
By (\ref{eq:mtdthreeblock}), a straightforward manipulation of conditional densities, 
and Fubini's theorem, we obtain 
\begin{eqnarray}
I 
&=& \int_{\mathbb{R}_{+}} \int_{\mathbb{R}_{+}} \int_{\mathbb{R}^p} 
\int_{\mathbb{R}^p} \int_{\mathbb{R}_{+}^p} \int_{\mathbb{R}_{+}^p} 
\pi\left(\tilde{\sigma}^2 \mid \tilde{\bm{\beta}},\bm{\tau},{\bf Y}\right)
\pi\left(\tilde{\bm{\beta}} \mid \bm{\tau},{\sigma^{2}},{\bf Y}\right)
\pi\left(\bm{\tau} \mid {\bm{\beta}},\sigma^{2},{\bf Y}\right) \nonumber\\
& & \pi\left(\bm{\beta} \mid \tilde{\bm{\tau}},{\sigma^{2}},{\bf Y}\right)
\pi\left({\sigma^{2}} \mid \tilde{\bm{\beta}},\tilde{\bm{\tau}},{\bf Y}\right)
\pi\left(\tilde{\bm{\tau}} \mid \tilde{\bm{\beta}},\tilde{\sigma}^2,{\bf Y}\right)
d\sigma^2 d\tilde{\sigma^2} d\bm{\beta} d\tilde{\bm{\beta}}
d{\bm{{\tau}}} d\tilde{\bm{\tau}} \label{eqtrace12}
\end{eqnarray}

\noindent
For convenience, we introduce and use the following notation in the subsequent proof. 
\begin{eqnarray}
& & \widehat{\bm{\beta}}={\bf A}_{\bm{\tau}}^{-1}X^T{\bf Y} \hspace{1.9in} ~\widehat{\bm{\beta}}_{*}={\bf A}_{\tilde{\bm{\tau}}}^{-1}X^T{\bf Y} \nonumber\\
& & {\Delta}_{1}=(\tilde{\bm{\beta}}-\widehat{\bm{\beta}})^T{\bf A}_{\bm{\tau}}(\tilde{\bm{\beta}}-\widehat{\bm{\beta}}) 
\hspace{1.2in} ~{\Delta}_{1*}=(\bm{\beta}-\widehat{\bm{\beta}}_{*})^T{\bf A}_{\tilde{\bm{\tau}}}(\bm{\beta} - 
\widehat{\bm{\beta}}_{*}) \nonumber\\
& & \tilde{\Delta}= ({\bf Y}-X\tilde{\bm{\beta}})^T({\bf Y}-X\tilde{\bm{\beta}})+{\tilde{\bm{\beta}}}^TD_{\bm{\tau}}^{-1}
{\tilde{\bm{\beta}}}+2\xi \hspace{0.2in} \tilde{\Delta}_{*}= ({\bf Y}-X\tilde{\bm{\beta}})^T({\bf Y}-X\tilde{\bm{\beta}})+
{\tilde{\bm{\beta}}}^TD_{\tilde{\bm{\tau}}}^{-1}{\tilde{\bm{\beta}}}+2\xi. \nonumber\\
& & \label{definition}
\end{eqnarray}

\noindent
We first show $I = \infty$ for the simpler case with $a > \frac{1}{2}$ and then consider 
the significantly more complicated case $0 < a \leq \frac{1}{2}$. 

\medskip

\noindent
{\bf Case 1: $a> 1/2$}

\noindent
Using $2K_{\nu}(x) \leq x^{-\nu} \Gamma(\nu) 2^\nu$ for $\nu > 0, x > 0$ 
(Proposition A7  of \cite{Pal:Khare:2014}), we obtain that if $a> \frac 1 2 ,$ 
\begin{equation} \label{eqtrace12.1}
\frac{\left(2b\sigma^2\right)^{\frac{a-\frac 1 2}{2}}}{2\abs*{\beta_j}^{a-\frac 1 2} 
K_{a-\frac 1 2}\left(\sqrt{2b\frac{\beta_j^2}{\sigma^2}}\right)} \geq 
\dfrac{(2b)^{a-\frac 1 2}}{\Gamma\big(a-\frac 1 2 \big) 2^{a-\frac 1 2}}. 
\end{equation}

\noindent
Similarly 
\begin{equation} \label{eqtrace12.2}
\frac{\left(2b\tilde{\sigma}^2\right)^{\frac{a-\frac 1 2}{2}}}{2\abs*{\tilde{\beta}_j}^{a-\frac 
1 2}K_{a-\frac 1 2}\left(\sqrt{2b\frac{\tilde{\beta}_j^2}{\tilde{\sigma}^2}}\right)}\geq
\dfrac{(2b)^{a-\frac 1 2}}{\Gamma\big(a-\frac 1 2 \big) 2^{a-\frac 1 2}}. 
\end{equation}

\noindent
Using (\ref{eq:cond1}), (\ref{eq:cond2}) and (\ref{eq:cond4}), along with 
(\ref{eqtrace12.1}) and (\ref{eqtrace12.2}), we get 
\begin{eqnarray}
& & ~\pi\left(\tilde{\sigma}^2 \mid \tilde{\bm{\beta}},\bm{\tau},{\bf Y}\right)
\pi\left(\tilde{\bm{\beta}} \mid \bm{\tau},{\sigma^{2}},{\bf Y}\right)
\pi\left(\bm{\tau} \mid {\bm{\beta}},\sigma^{2},{\bf Y}\right) \nonumber\\
& & ~\pi\left(\bm{\beta} \mid \tilde{\bm{\tau}},{\sigma^{2}},{\bf Y}\right) 
\pi\left({\sigma^{2}} \mid \tilde{\bm{\beta}},\tilde{\bm{\tau}},{\bf Y}\right)
\pi\left(\tilde{\bm{\tau}} \mid {\tilde{\bm{\beta}}},\tilde{\sigma}^2,{\bf Y}\right) 
\nonumber\\
&=& D_1\Bigg\{\frac{{\tilde{\Delta}}^{\frac{n+p+2\alpha}{2}} \exp \left(-\frac{{\tilde{\Delta}}}{2\tilde{\sigma}^2} \right)}
{{(\tilde{{\sigma}}^2)}^{{\frac{n+p+2\alpha}{2}}+1}} \Bigg\}\Bigg\{
\frac{{\abs*{{{\bf A}_{\bm{\tau}}}}}^{\frac{1}{2}} \exp \left(-\frac{{{\Delta}}_1}{2\sigma^2} \right)}{\sigma^p}
\Bigg\} \nonumber\\
& &\Bigg\{
\prod\limits_{j=1}^p\frac{\left(2b\sigma^2\right)^{\frac{a-\frac 1 2}{2}}}{2\abs*{\beta_j}^{a-\frac 1 2}K_{a-\frac 1 2}\left(\sqrt{2b\frac{\beta_j^2}{\sigma^2}}\right)}\tau_j^{\left(a-\frac 1 2\right)-1}e^{-\frac 1 2 \left\{2b\tau_j+\frac{\beta_j^2}{\sigma^2}\frac{1}{\tau_j}\right\} }
\Bigg\} 
\Bigg\{
\frac{{\abs*{{{\bf A}_{\bm{\tilde{\tau}}}}}}^{\frac 1 2 } \exp \left(-\frac{{{\Delta}}_{1*}}{2\sigma^2} \right)}{\sigma^p}
\Bigg\} \nonumber\\
& &\Bigg\{\frac{{\tilde{\Delta}}_{*}^{\frac{n+p+2\alpha}{2}} \exp \left(-\frac{{\tilde{\Delta}}_{*}}{2{\sigma^2}} \right)}
{{({\sigma^2})}^{{\frac{n+p+2\alpha}{2}}+1}}\Bigg\} 
\Bigg\{
\prod\limits_{j=1}^p\frac{\left(2b\tilde{\sigma}^2\right)^{\frac{a-\frac 1 2}{2}}}{2\abs*{\tilde{\beta}_j}^{a-\frac 1 2}K_{a-\frac 1 2}\left(\sqrt{2b\frac{\tilde{\beta}_j^2}{\tilde{\sigma}^2}}\right)}\tilde{\tau}_j^{\left(a-\frac 1 2\right)-1}e^{-\frac 1 2 \left\{2b\tilde{\tau}_j+\frac{\tilde{\beta}_j^2}{\tilde{\sigma}^2}\frac{1}{\tilde{\tau}_j}\right\} }\Bigg\} \nonumber\\
&\geq& D_1 f_1(\bm{\tau},\tilde{\bm{\tau}})\Bigg\{\frac{{\tilde{\Delta}}^{\frac{n+p+2\alpha}{2}} \exp \left(-\frac{{\tilde{\Delta}} + 
		{\tilde{\bm{\beta}}}^TD_{\tilde{\bm{\tau}}}^{-1}\tilde{\bm{\beta}}}{2\tilde{{\sigma}}^2} \right)}{{(\tilde{\sigma}^2)}^{{\frac{n+p+2\alpha}{2}}+1}}
\Bigg\}\Bigg\{\frac{{\tilde{\Delta}}_{*}^{\frac{n+p+2\alpha}{2}} \exp \left(-\frac{{{\Delta}_{1}+\Delta_{1*}+\tilde{\Delta}_{*} + 
			\bm{\beta}^TD_{\bm{\tau}}^{-1}\bm{\beta}}}{2{\sigma^2}} \right)}{{({\sigma^2})}^{{\frac{n+p+2\alpha}{2}}+p+1}}\Bigg\}\label{eqtrace13}
\end{eqnarray}

\noindent
where 
$$
D_1=\frac{1}{\left[(2\pi)^{{p}}2^{\frac{n+p+2\alpha}{2}}\Gamma(\frac{n+p+2\alpha}{2})\right]^2}\left(	\dfrac{(2b)^{a-\frac 1 2}}{\Gamma\big(a-\frac 1 2 \big) 2^{a-\frac 1 2}}\right)^{2p} 
$$

\noindent
and
$$
f_1(\bm{\tau},\tilde{\bm{\tau}})=\Bigg\{
\prod\limits_{j=1}^{p} ~\tau_j^{\left(a-\frac 1 2\right)-1}e^{-b\tau_j}
\Bigg\}\Bigg\{
\prod\limits_{j=1}^{p} ~\tilde{\tau}_j^{\left(a-\frac 1 2\right)-1}e^{-b\tilde{\tau}_j}
\Bigg\}{\abs*{{{\bf A}_{\bm{\tau}}}}}^{\frac 1 2 }{\abs*{{{\bf A}_{\bm{\tilde{\tau}}}}}}^{\frac 1 2 }. 
$$

\noindent
It follows from (\ref{eqtrace13}) that
\begin{eqnarray*}
	& & \int_{\mathbb{R}_{+}}\int_{\mathbb{R}_{+}}\int_{\mathbb{R}^p}\int_{\mathbb{R}^p} \pi\left(\tilde{\sigma}^2 \mid 
	\tilde{\bm{\beta}},\bm{\tau},{\bf Y}\right) \pi\left(\tilde{\bm{\beta}} \mid \bm{\tau},{\sigma^{2}},{\bf Y}\right) 
	\pi\left(\bm{\tau} \mid {\bm{\beta}},\sigma^{2},{\bf Y}\right)\\
	& & \pi\left(\bm{\beta} \mid \tilde{\bm{\tau}},{\sigma^{2}},{\bf Y}\right) 
	\pi \left({\sigma^{2}} \mid \tilde{\bm{\beta}},\tilde{\bm{\tau}}, 
	{\bf Y}\right) \pi\left(\tilde{\bm{\tau}} \mid \tilde{\bm{\beta}},\tilde{\sigma}^2,{\bf Y}\right)d\sigma^2d \tilde{\sigma}^2d\bm{\beta} 
	d\tilde{\bm{\beta}}\\
	&\geq&
	D_1 f_1(\bm{\tau},\tilde{\bm{\tau}})
	\int_{\mathbb{R}_{+}}\int_{\mathbb{R}_{+}}\int_{\mathbb{R}^p}\int_{\mathbb{R}^p}
	\Bigg\{\frac{{\tilde{\Delta}}^{\frac{n+p+2\alpha}{2}} \exp \left(-\frac{{\tilde{\Delta}} + 
			{\tilde{\bm{\beta}}}^TD_{\tilde{\bm{\tau}}}^{-1}\tilde{\bm{\beta}}}{2\tilde{{\sigma}}^2} \right)}{{(\tilde{\sigma}^2)}^{{\frac{n+p+2\alpha}{2}}+1}}
	\Bigg\}\Bigg\{\frac{{\tilde{\Delta}}_{*}^{\frac{n+p+2\alpha}{2}} \exp \left(-\frac{{{\Delta}_{1}+\Delta_{1*}+\tilde{\Delta}_{*} + 
				\bm{\beta}^TD_{\bm{\tau}}^{-1}\bm{\beta}}}{2{\sigma^2}} \right)}{{({\sigma^2})}^{{\frac{n+p+2\alpha}{2}}+p+1}}\Bigg\}\\
	& & d\sigma^2d \tilde{\sigma}^2d\bm{\beta} 
	d\tilde{\bm{\beta}}\\
	&\stackrel{(a)}{=}& \infty 
\end{eqnarray*}
\noindent
for every $\left(\bm{\tau},\tilde{\bm{\tau}}\right) \in {\mathbb{R}_{+}^{p}} \times 
{\mathbb{R}_{+}^{p}}$. Here (a) follows by repeating verbatim the arguments between 
Equations (S4) - (S12) in \cite{SRKZ:2016}.  We conclude from this fact that the 
Markov operator $\tilde{K}$ is not Hilbert-Schmidt when $a > \frac{1}{2}$. 

\medskip

\noindent
{\bf Case 2: $0 <a \leq 1/2$}

\noindent
By the integral formula (see \cite{Abramowitz:Stegun:1965}, Page 376)
\[K_{\nu}(t)= \int_{0}^{\infty} e^{-t\cosh z} \cosh(vz) \,dz, \nu\in\mathcal{R}. 
\]

\noindent
Since the hyperbolic function $\cosh$ is strictly decreasing on interval $(-\infty,0]$, for 
every $x > 0$, $K_{\nu} (x)$ is strictly decreasing as ${\nu}$ increases on the interval 
$(-\infty,0].$  Note that when $0<a\leq \frac 1 2,$ $-a-\frac 3 2 < {a-\frac 1 2} \leq 0$. It 
follows that 
$$
K_{a-\frac 1 2}(x) < K_{-a-\frac 3 2}(x) 
$$

\noindent
for all $x>0$. Moreover, when $\nu<0 ~\text{and~} x>0, $ $2K_{\nu}(x) \leq x^{{\nu}} 
\Gamma(-{\nu}) 2^{-{\nu}}$ (see Proposition A7 of \cite{Pal:Khare:2014}), which 
implies 
\[2K_{a-\frac 1 2}\left(\sqrt{2b\frac{\beta_j^2}{\sigma^2}}\right) <   2K_{-a-\frac 3 
2}\left(\sqrt{2b\frac{\beta_j^2}{\sigma^2}}\right) \leq  \left(\sqrt{2b\frac{\beta_j^2}
{\sigma^2}}\right)^{-a-\frac 3 2}\Gamma({a+\frac 3 2}) ~2^{{a+\frac 3 2}} 
\]

\noindent
and 
\begin{equation} \label{eqtrace13.1}
\frac{\left(2b\sigma^2\right)^{\frac{a-\frac 1 2}{2}}}{2\abs*{\beta_j}^{a-\frac 1 2}K_{a-
\frac 1 2}\left(\sqrt{2b\frac{\beta_j^2}{\sigma^2}}\right)} \geq
\frac{b^{{a+\frac 1 2}}}{2\Gamma({a+\frac 3 2})}
\left(\frac{\beta_j^2}{\sigma^2}\right). 
\end{equation}

\noindent
Similarly, we get
\begin{equation} \label{eqtrace13.2}
\frac{\left(2b\tilde{\sigma}^2\right)^{\frac{a-\frac 1 2}{2}}}{2\abs*{\tilde{\beta}_j}^{a-\frac 
1 2}K_{a-\frac 1 2}\left(\sqrt{2b\frac{\tilde{\beta}_j^2}{\tilde{\sigma}^2}}\right)}
\geq \frac{b^{{a+\frac 1 2}}}{2\Gamma({a+\frac 3 2})}
\left(\frac{\tilde{\beta}_j^2}{\tilde{\sigma}^2}\right). 
\end{equation}

\noindent
Using (\ref{eq:cond1}), (\ref{eq:cond2}) and (\ref{eq:cond4}), along with 
(\ref{eqtrace13.1}) and (\ref{eqtrace13.2}), we obtain 
\begin{eqnarray}
& &~ \pi \left(\tilde{\sigma}^2 \mid \tilde{\bm{\beta}},\bm{\tau},{\bf Y}\right)
\pi \left(\tilde{\bm{\beta}} \mid \bm{\tau},{\sigma^{2}},{\bf Y}\right)
\pi \left(\bm{\tau} \mid {\bm{\beta}},\sigma^{2},{\bf Y}\right)
\nonumber\\
& & ~\pi \left(\bm{\beta} \mid \tilde{\bm{\tau}},{\sigma^{2}},{\bf Y}\right) 
\pi \left({\sigma^{2}} \mid \tilde{\bm{\beta}},\tilde{\bm{\tau}},{\bf Y}\right)
\pi \left(\tilde{\bm{\tau}} \mid {\tilde{\bm{\beta}}},\tilde{\sigma}^2,{\bf Y}\right) 
\nonumber\\
&=& D_2\Bigg\{\frac{{\tilde{\Delta}}^{\frac{n+p+2\alpha}{2}} \exp \left(-\frac{{\tilde{\Delta}}}{2\tilde{\sigma}^2} \right)}
{{(\tilde{{\sigma}}^2)}^{{\frac{n+p+2\alpha}{2}}+1}} \Bigg\}\Bigg\{
\frac{{\abs*{{{\bf A}_{\bm{\tau}}}}}^{\frac{1}{2}} \exp \left(-\frac{{{\Delta}}_1}{2\sigma^2} \right)}{\sigma^p}
\Bigg\} \nonumber\\
& &\Bigg\{
\prod\limits_{j=1}^p\frac{\left(2b\sigma^2\right)^{\frac{a-\frac 1 2}{2}}}{2\abs*{\beta_j}^{a-\frac 1 2}K_{a-\frac 1 2}\left(\sqrt{2b\frac{\beta_j^2}{\sigma^2}}\right)}\tau_j^{\left(a-\frac 1 2\right)-1}e^{-\frac 1 2 \left\{2b\tau_j+\frac{\beta_j^2}{\sigma^2}\frac{1}{\tau_j}\right\} }
\Bigg\} 
\Bigg\{
\frac{{\abs*{{{\bf A}_{\bm{\tilde{\tau}}}}}}^{\frac 1 2 } \exp \left(-\frac{{{\Delta}}_{1*}}{2\sigma^2} \right)}{\sigma^p}
\Bigg\} \nonumber \\
& &\Bigg\{\frac{{\tilde{\Delta}}_{*}^{\frac{n+p+2\alpha}{2}} \exp \left(-\frac{{\tilde{\Delta}}_{*}}{2{\sigma^2}} \right)}
{{({\sigma^2})}^{{\frac{n+p+2\alpha}{2}}+1}}\Bigg\} 
\Bigg\{
\prod\limits_{j=1}^p\frac{\left(2b\tilde{\sigma}^2\right)^{\frac{a-\frac 1 2}{2}}}{2\abs*{\tilde{\beta}_j}^{a-\frac 1 2}K_{a-\frac 1 2}\left(\sqrt{2b\frac{\tilde{\beta}_j^2}{\tilde{\sigma}^2}}\right)}\tilde{\tau}_j^{\left(a-\frac 1 2\right)-1}e^{-\frac 1 2 \left\{2b\tilde{\tau}_j+\frac{\tilde{\beta}_j^2}{\tilde{\sigma}^2}\frac{1}{\tilde{\tau}_j}\right\} }\Bigg\} \nonumber\\
& \geq& D_2 f_1(\bm{\tau},\tilde{\bm{\tau}})
\prod\limits_{j=1}^{p}\left({\beta_j^2\tilde{\beta}_j^2}\right)
\Bigg\{\frac{ \exp \left(-\frac{{\tilde{\Delta}} + 
		{\tilde{\bm{\beta}}}^TD_{\tilde{\bm{\tau}}}^{-1}\tilde{\bm{\beta}}}{2\tilde{{\sigma}}^2} \right)}{{(\tilde{\sigma}^2)}^{\frac{n+p+2\alpha}{2}+p+1}}
\Bigg\}\Bigg\{\frac{ \exp \left(-\frac{{{\Delta}_{1}+\Delta_{1*}+\tilde{\Delta}_{*} + 
			\bm{\beta}^TD_{\bm{\tau}}^{-1}\bm{\beta}}}{2{\sigma^2}} \right)}{{({\sigma^2})}^{\frac{n+p+2\alpha}{2}+2p+1}}\Bigg\},\label{eqtrace14}
\end{eqnarray}

\noindent
where 
$$
D_2=\frac{(2\xi)^{n+p+2\alpha}}{\left[(2\pi)^{{p}}2^{\frac{n+p+2\alpha}{2}}
\Gamma(\frac{n+p+2\alpha}{2})\right]^2} 
\left( \frac{b^{{a+\frac 1 2}}}{2\Gamma({a+\frac 3 2})} \right)^{2p}, 
$$

$$
f_1(\bm{\tau},\tilde{\bm{\tau}})=\Bigg\{ \prod\limits_{j=1}^{p} 
~\tau_j^{a-2}e^{-b\tau_j} \Bigg\} \Bigg\{ 
\prod\limits_{j=1}^{p} ~\tilde{\tau}_j^{a-2} e^{-b\tilde{\tau}_j}
\Bigg\}, 
$$

\noindent
and the last inequality follows by 
\[{\tilde{\Delta}}^{\frac{n+p+2\alpha}{2}} \geq (2\xi)^{\frac{n+p+2\alpha}
{2}},~{\tilde{\Delta}}_{*}^{\frac{n+p+2\alpha}{2}}\geq (2\xi)^{\frac{n+p+2\alpha}{2}},~
{\abs*{{{\bf A}_{\bm{\tau}}}}}^{\frac 1 2 }\geq |D_{\bm{\tau}}|^{-\frac 1 2}~\text{and}~
{\abs*{{{\bf A}_{\bm{\tilde{\tau}}}}}}^{\frac 1 2 }\geq |D_{\bm{\tilde{\tau}}}|^{-\frac 1 2}. 
\]

\noindent
It follows by (\ref{eqtrace14}) and the form of the Inverse-Gamma density that
\begin{eqnarray}
& & ~\int_{\mathbb{R}_{+}}\int_{\mathbb{R}_{+}} \pi \left(\tilde{\sigma}^2~|~
\tilde{\bm{\beta}},\bm{\tau},{\bf Y}\right) \pi \left(\tilde{\bm{\beta}}~|~\bm{\tau},
{\sigma^{2}},{\bf Y}\right) \pi \left(\bm{\tau}~|~{\bm{\beta}},\sigma^{2},{\bf Y}\right)  
\nonumber\\
& & \pi \left(\bm{\beta}~|~ \tilde{\bm{\tau}},{\sigma^{2}},{\bf Y}\right)
\pi \left({\sigma^{2}}~|~\tilde{\bm{\beta}},\tilde{\bm{\tau}},{\bf Y}\right)
\pi \left(\tilde{\bm{\tau}}~|~\tilde{\bm{\beta}},\tilde{\sigma}^2,{\bf Y}\right)d\sigma^2d
\tilde{{\sigma^2}} \nonumber\\
&\geq& D_3 f_1(\bm{\tau},\tilde{\bm{\tau}})
\prod\limits_{j=1}^{p}\left({\beta_j^2\tilde{\beta}_j^2}\right) 
\left\{\frac{1}{\Big[{{\tilde{\Delta}} + 
		{\tilde{\bm{\beta}}}^TD_{\tilde{\bm{\tau}}}^{-1}\tilde{\bm{\beta}}}\Big]^{{\frac{n+p+2\alpha}{2}+p}}}\right\}
\left\{\frac{1}{\Big [{{{{\Delta}_{1}+\Delta_{1*}+\tilde{\Delta}_{*}+\bm{\beta}^T 
				D_{\bm{\tau}}^{-1}\bm{\beta}}}}\Big ]^{\frac{n+p+2\alpha}{2}+2p}}\right\} \nonumber\\ & \label{eqtrace15}
\end{eqnarray}

\noindent
where 
$$
D_3=2^{n+4p+2\alpha}~\Gamma\left(\frac{n+p+2\alpha}{2}+p\right)~\Gamma\left(\frac{n+p+2\alpha}{2}+2p\right)~D_2
$$

\noindent
We now establish some inequalities which will help converting the lower bound in 
(\ref{eqtrace15}) into a simpler form. By (\ref{definition}), it follows that 
\begin{align*}
{{\tilde{\Delta}} + 
	{\tilde{\bm{\beta}}}^TD_{\tilde{\bm{\tau}}}^{-1}\tilde{\bm{\beta}}}
& = \tilde{\bm{\beta}}^T(X^TX+D_{\bm{\tau}}^{-1})\tilde{\bm{\beta}}-2{\bf Y}^TX\tilde{\bm{\beta}}
+{\bf Y}^T{\bf Y}+2\xi+	{\tilde{\bm{\beta}}}^TD_{\tilde{\bm{\tau}}}^{-1}\tilde{\bm{\beta}}\nonumber\\ & 
\leq  \tilde{\bm{\beta}}^T(X^TX+D_{\bm{\tau}}^{-1})\tilde{\bm{\beta}}-2{\bf Y}^TX\tilde{\bm{\beta}}
+{\bf Y}^T{\bf Y} + \tilde{\Delta}_{*}\nonumber\\
& \leq \tilde{\bm{\beta}}^T(X^TX+D_{\bm{\tau}}^{-1})\tilde{\bm{\beta}}-2{\bf Y}^TX\tilde{\bm{\beta}}
+{\bf Y}^T{\bf Y} + \tilde{\Delta}_{*}+\Delta_{1*}+\bm{\beta}^T 
D_{\bm{\tau}}^{-1}\bm{\beta}, 
\end{align*}

\noindent
and
\begin{align*}
& ~~~~{{{{\Delta}_{1}+\Delta_{1*}+\tilde{\Delta}_{*}+\bm{\beta}^T 
			D_{\bm{\tau}}^{-1}\bm{\beta}}}} \nonumber \\ 
& =
\tilde{\bm{\beta}}^T(X^TX+D_{\bm{\tau}}^{-1})\tilde{\bm{\beta}}-2{\bf Y}^TX\tilde{\bm{\beta}}
+{\bf Y}^TX(X^TX+D_{\bm{\tau}}^{-1})^{-1}X^T{\bf Y}
+\Delta_{1*}+\tilde{\Delta}_{*}+\bm{\beta}^T 
D_{\bm{\tau}}^{-1}\bm{\beta} \nonumber \\
& \leq  \tilde{\bm{\beta}}^T(X^TX+D_{\bm{\tau}}^{-1})\tilde{\bm{\beta}}-2{\bf Y}^TX\tilde{\bm{\beta}}
+{\bf Y}^T{\bf Y} +  \tilde{\Delta}_{*}+\Delta_{1*}+\bm{\beta}^T 
D_{\bm{\tau}}^{-1}\bm{\beta}. 
\end{align*}

\noindent
Also, note that 
\begin{align*}
& ~~~~\tilde{\bm{\beta}}^T(X^TX+D_{\bm{\tau}}^{-1})\tilde{\bm{\beta}}-2{\bf Y}^TX\tilde{\bm{\beta}}
+{\bf Y}^T{\bf Y} +  \tilde{\Delta}_{*}+\Delta_{1*}+\bm{\beta}^T 
D_{\bm{\tau}}^{-1}\bm{\beta}  \nonumber \\ 
& = \tilde{\bm{\beta}}^T(X^TX+D_{\bm{\tau}}^{-1})\tilde{\bm{\beta}}-2{\bf Y}^TX\tilde{\bm{\beta}}
+{\bf Y}^T{\bf Y} \nonumber\\
& ~~~+ \tilde{\bm{\beta}}^T(X^TX+D_{\tilde{\bm{\tau}}}^{-1})\tilde{\bm{\beta}}-2{\bf Y}^TX\tilde{\bm{\beta}}+{\bf Y}^T{\bf Y}+2\xi+\Delta_{1*} +\bm{\beta}^T 
D_{\bm{\tau}}^{-1}\bm{\beta}  \nonumber  \\
& = \tilde{\bm{\beta}}^T(2X^TX+D_{\bm{\tau}}^{-1}+D_{\tilde{\bm{\tau}}}^{-1})\tilde{\bm{\beta}}-4{\bf Y}^TX\tilde{\bm{\beta}}
+2{\bf Y}^T{\bf Y}+2\xi+ {\Delta}_{1*}+\bm{\beta}^TD_{\bm{\tau}}^{-1}\bm{\beta}  \nonumber  \\
& = (\tilde{\bm{\beta}}-\mu)^T(2X^TX+D_{\bm{\tau}}^{-1}+D_{\tilde{\bm{\tau}}}^{-1})(\tilde{\bm{\beta}}-\mu)-4{\bf Y}^TX(2X^TX+D_{\bm{\tau}}^{-1}+D_{\tilde{\bm{\tau}}}^{-1})^{-1}X^T{\bf Y}\nonumber\\
&~~~ +f_2(\bm{\beta},\tilde{\bm{\tau}})+\bm{\beta}^TD_{\bm{\tau}}^{-1}\bm{\beta}\nonumber\\
& \leq (\tilde{\bm{\beta}}-\mu)^T(2X^TX+D_{\bm{\tau}}^{-1}+D_{\tilde{\bm{\tau}}}^{-1})(\tilde{\bm{\beta}}-\mu)+f_2(\bm{\beta},\tilde{\bm{\tau}})+\bm{\beta}^TD_{\bm{\tau}}^{-1}\bm{\beta}\nonumber\\
& \leq \left[(\tilde{\bm{\beta}}-\mu)^T(2X^TX+D_{\bm{\tau}}^{-1}+D_{\tilde{\bm{\tau}}}^{-1})(\tilde{\bm{\beta}}-\mu)+f_2(\bm{\beta},\tilde{\bm{\tau}})+1\right](\bm{\beta}^TD_{\bm{\tau}}^{-1}\bm{\beta}+1)
\end{align*}

\noindent
where 
$$
\mu=(2X^TX+D_{\bm{\tau}}^{-1}+D_{\tilde{\bm{\tau}}}^{-1})^{-1}X^T {\bf Y} 
\mbox{ and } f_2(\bm{\beta},\tilde{\bm{\tau}})= 2{\bf Y}^T{\bf Y}+2\xi+\Delta_{1*}. 
$$

\noindent
By (\ref{eqtrace15}), we get 
\begin{eqnarray*}
	& & ~\int_{\mathbb{R}_{+}}\int_{\mathbb{R}_{+}} \pi \left(\tilde{\sigma}^2~|~\tilde{\bm{\beta}},\bm{\tau},{\bf Y}\right)
	\pi \left(\tilde{\bm{\beta}}~|~\bm{\tau},{\sigma^{2}},{\bf Y}\right)
	\pi \left(\bm{\tau}~|~{\bm{\beta}},\sigma^{2},{\bf Y}\right)  \nonumber\\
	& & \pi \left(\bm{\beta}~|~ \tilde{\bm{\tau}},{\sigma^{2}},{\bf Y}\right)
	\pi \left({\sigma^{2}}~|~\tilde{\bm{\beta}},\tilde{\bm{\tau}},{\bf Y}\right)
	\pi \left(\tilde{\bm{\tau}}~|~\tilde{\bm{\beta}},\tilde{\sigma}^2,{\bf Y}\right)d\sigma^2d\tilde{{\sigma^2}} \nonumber\\
	&\geq& \frac{D_2f_1(\bm{\tau},\tilde{\bm{\tau}})}{(\bm{\beta}^TD_{\bm{\tau}}^{-1}\bm{\beta}+1)^{n+4p+2\alpha}}
	\frac{ \prod\limits_{j=1}^{p}\left({\beta_j^2\tilde{\beta}_j^2}\right) }{\left[(\tilde{\bm{\beta}}-\mu)^T(2X^TX+D_{\bm{\tau}}^{-1}+D_{\tilde{\bm{\tau}}}^{-1})(\tilde{\bm{\beta}}-\mu)+f_2(\bm{\beta},\tilde{\bm{\tau}})+1\right]^{n+4p+2\alpha}}\nonumber\\
	& =& \frac{D_2f_1(\bm{\tau},\tilde{\bm{\tau}})}{(\bm{\beta}^TD_{\bm{\tau}}^{-1}\bm{\beta}+1)^{n+4p+2\alpha}}
	\frac{ \prod\limits_{j=1}^{p}\left({\beta_j^2\tilde{\beta}_j^2}\right) }{\left[\left(f_2(\bm{\beta},\tilde{\bm{\tau}})+1\right)\left(1+\frac{(\tilde{\bm{\beta}}-\mu)^T(2X^TX+D_{\bm{\tau}}^{-1}+D_{\tilde{\bm{\tau}}}^{-1})(\tilde{\bm{\beta}}-\mu)}{f_2(\bm{\beta},\tilde{\bm{\tau}})+1}\right)\right]^{n+4p+2\alpha}}\nonumber\\
	& =& \frac{D_2f_1(\bm{\tau},\tilde{\bm{\tau}})}{\left(f_2(\bm{\beta},\tilde{\bm{\tau}})+1\right)^{n+4p+2\alpha}(\bm{\beta}^TD_{\bm{\tau}}^{-1}\bm{\beta}+1)^{n+4p+2\alpha}}
	\frac{ \prod\limits_{j=1}^{p}\left({\beta_j^2\tilde{\beta}_j^2}\right) }{\left(1+\frac{(\tilde{\bm{\beta}}-\mu)^T(2X^TX+D_{\bm{\tau}}^{-1}+D_{\tilde{\bm{\tau}}}^{-1})(\tilde{\bm{\beta}}-\mu)}{f_2(\bm{\beta},\tilde{\bm{\tau}})+1}\right)^{n+4p+2\alpha}}\nonumber\\
	& \geq & \frac{D_2f_1(\bm{\tau},\tilde{\bm{\tau}})}{\left(f_2(\bm{\beta},\tilde{\bm{\tau}})+1\right)^{n+4p+2\alpha}(\bm{\beta}^TD_{\bm{\tau}}^{-1}\bm{\beta}+1)^{n+4p+2\alpha}}
	\frac{ \prod\limits_{j=1}^{p}\left({\beta_j^2\tilde{\beta}_j^2}\right) }{\left(1+\frac{(\tilde{\bm{\beta}}-\mu)^T(2\lambda {I}_p+D_{\bm{\tau}}^{-1}+D_{\tilde{\bm{\tau}}}^{-1})(\tilde{\bm{\beta}}-\mu)}{f_2(\bm{\beta},\tilde{\bm{\tau}})+1}\right)^{n+4p+2\alpha}}\nonumber\\
\end{eqnarray*}

\begin{eqnarray}
& \geq & \frac{D_2f_1(\bm{\tau},\tilde{\bm{\tau}})}{\left(f_2(\bm{\beta},\tilde{\bm{\tau}})+1\right)^{n+4p+2\alpha}(\bm{\beta}^TD_{\bm{\tau}}^{-1}\bm{\beta}+1)^{n+4p+2\alpha}}
\frac{ \prod\limits_{j=1}^{p}\left({\beta_j^2\tilde{\beta}_j^2}\right) }{\prod\limits_{j=1}^{p}\left(1+\frac{\left(\tilde{{\beta}}_j-\mu_j\right)^2}{\nu_j\epsilon_j^2}\right)^{n+4p+2\alpha}}\nonumber\\
& = & \frac{D_2f_1(\bm{\tau},\tilde{\bm{\tau}})\prod\limits_{j=1}^{p}{\beta_j^2\prod\limits_{j=2}^{p}\tilde{\beta}_j^2}}{\left(f_2(\bm{\beta},\tilde{\bm{\tau}})+1\right)^{n+4p+2\alpha}{\prod\limits_{j=2}^{p}\left(1+\frac{\left(\tilde{{\beta}}_j-\mu_j\right)^2}{\nu_j\epsilon_j^2}\right)^{n+4p+2\alpha}}(\bm{\beta}^TD_{\bm{\tau}}^{-1}\bm{\beta}+1)^{n+4p+2\alpha}}~
\frac{\tilde{\beta}_1^2} {\left(1+\frac{\left(\tilde{{\beta}}_1-\mu_1\right)^2}{\nu_1\epsilon_1^2}\right)^{\frac{1+\nu_1}{2}}}\nonumber\\\label{eqtrace16}
\end{eqnarray}
where $\mu_j=e_j^T\mu,~\nu_j=2n+8p+4\alpha-1,~\epsilon_j=\sqrt{\frac{f_2(\bm{\beta},\tilde{\bm{\tau}})+1}{(2\lambda+\frac{1}{\tau_j}+\frac{1}{\tilde{\tau}_j})\nu_j}}$ and $\lambda$ is the greatest eigenvalue of matrix $X^TX.$ By Proposition A4, we have
\begin{eqnarray}
\int_{\mathbb{R}}\frac{\tilde{\beta}_1^2} {\left(1+\frac{\left(\tilde{\bm{\beta}}_1-\mu_1\right)^2}{\nu_1\epsilon_1^2}\right)^{\frac{1+\nu_1}{2}}}d{\tilde{{\beta}}_1}\geq f_3(\bm{\beta},\tilde{\bm{\tau}}) \left({2\lambda+\frac{1}{\tau_1}+\frac{1}{\tilde{\tau}_1}}\right)^{-\frac 3 2}\label{eqtrace17}
\end{eqnarray}
\noindent Hence, it follows from (\ref{eqtrace16}) and (\ref{eqtrace17}) that 
\begin{eqnarray}
& & ~~\int_{\mathbb{R}_{+}}\int_{\mathbb{R}_{+}}\int_{\mathbb{R}} f\left(\tilde{\sigma}^2~|~\tilde{\bm{\beta}},\bm{\tau},{\bf Y}\right)
f\left(\tilde{\bm{\beta}}~|~\bm{\tau},{\sigma^{2}},{\bf Y}\right)
f\left(\bm{\tau}~|~{\bm{\beta}},\sigma^{2},{\bf Y}\right)  \nonumber\\
& &~~ f\left(\bm{\beta}~|~ \tilde{\bm{\tau}},{\sigma^{2}},{\bf Y}\right)
f\left({\sigma^{2}}~|~\tilde{\bm{\beta}},\tilde{\bm{\tau}},{\bf Y}\right)
f\left(\tilde{\bm{\tau}}~|~\tilde{\bm{\beta}},\tilde{\sigma}^2,{\bf Y}\right)d\sigma^2d\tilde{{\sigma^2}}d\tilde{{\beta}}_1 \nonumber\\
& & \geq 
\frac{D_2~f_1(\bm{\tau},\tilde{\bm{\tau}})\prod\limits_{j=1}^{p}{\beta_j^2\prod\limits_{j=2}^{p}\tilde{\beta}_j^2}}{\left(f_2(\bm{\beta},\tilde{\bm{\tau}})+1\right)^{n+4p+2\alpha}{\prod\limits_{j=2}^{p}\left(1+\frac{\left(\tilde{{\beta}}_j-\mu_j\right)^2}{\nu_j\epsilon_j^2}\right)^{n+4p+2\alpha}}(\bm{\beta}^TD_{\bm{\tau}}^{-1}\bm{\beta}+1)^{n+4p+2\alpha}}~\nonumber\\
& & ~~~~\times  \int_{\mathbb{R}}\frac{\tilde{\beta}_1^2} {\left(1+\frac{\left(\tilde{\bm{\beta}}_1-\mu_1\right)^2}{\nu_1\epsilon_1^2}\right)^{\frac{1+\nu_1}{2}}}~d{\tilde{{\beta}}_1}\nonumber\\
& & \geq 
\frac{D_2~f_1(\bm{\tau},\tilde{\bm{\tau}})\prod\limits_{j=1}^{p}{\beta_j^2\prod\limits_{j=2}^{p}\tilde{\beta}_j^2}~f_3(\bm{\beta},\tilde{\bm{\tau}}) \left({2\lambda+\frac{1}{\tau_1}+\frac{1}{\tilde{\tau}_1}}\right)^{-\frac 3 2}}{\left(f_2(\bm{\beta},\tilde{\bm{\tau}})+1\right)^{n+4p+2\alpha}{\prod\limits_{j=2}^{p}\left(1+\frac{\left(\tilde{{\beta}}_j-\mu_j\right)^2}{\nu_j\epsilon_j^2}\right)^{n+4p+2\alpha}}(\bm{\beta}^TD_{\bm{\tau}}^{-1}\bm{\beta}+1)^{n+4p+2\alpha}}\nonumber\\
& & =
\frac{D_2~\Bigg\{
	\prod\limits_{j=2}^{p} ~\tau_j^{\left(a-\frac 1 2\right)-1}e^{-b\tau_j}
	\Bigg\}\Bigg\{
	\prod\limits_{j=1}^{p} ~\tilde{\tau}_j^{\left(a-\frac 1 2\right)-1}e^{-b\tilde{\tau}_j}
	\Bigg\}\prod\limits_{j=1}^{p}{\beta_j^2\prod\limits_{j=2}^{p}\tilde{\beta}_j^2}~f_3(\bm{\beta},\tilde{\bm{\tau}}) }{\left(f_2(\bm{\beta},\tilde{\bm{\tau}})+1\right)^{n+4p+2\alpha}{\prod\limits_{j=2}^{p}\left(1+\frac{\left(\tilde{{\beta}}_j-\mu_j\right)^2}{\nu_j\epsilon_j^2}\right)^{n+4p+2\alpha}}}~\nonumber\\
&&~~~~\times~~\frac{
	\tau_1^{\left(a-\frac 1 2\right)-1}e^{-b\tau_1}
	\left({2\lambda+\frac{1}{\tau_1}+\frac{1}{\tilde{\tau}_1}}\right)^{-\frac 3 2}}{
	\left(   	\sum\limits_{j=2}^{p}\frac{\beta_j^2}{\tau_j}+\frac{\beta_1^2}{\tau_1}+1\right)^{n+4p+2\alpha}}\label{eqtrace18}
\end{eqnarray}
From (\ref{eqtrace18}), we obtain 
\begin{eqnarray}
& & ~~\int_{\mathbb{R}_{+}}\int_{\mathbb{R}_{+}}\int_{\mathbb{R}} \int_{\mathbb{R}_{+}} f\left(\tilde{\sigma}^2~|~\tilde{\bm{\beta}},\bm{\tau},{\bf Y}\right)
f\left(\tilde{\bm{\beta}}~|~\bm{\tau},{\sigma^{2}},{\bf Y}\right)
f\left(\bm{\tau}~|~{\bm{\beta}},\sigma^{2},{\bf Y}\right)  \nonumber\\
& &~~ f\left(\bm{\beta}~|~ \tilde{\bm{\tau}},{\sigma^{2}},{\bf Y}\right)
f\left({\sigma^{2}}~|~\tilde{\bm{\beta}},\tilde{\bm{\tau}},{\bf Y}\right)
f\left(\tilde{\bm{\tau}}~|~\tilde{\bm{\beta}},\tilde{\sigma}^2,{\bf Y}\right)d\sigma^2d\tilde{{\sigma^2}}d\tilde{{\beta}}_1d\tau_1 \nonumber\\
& & \geq \frac{D_2~\Bigg\{
	\prod\limits_{j=2}^{p} ~\tau_j^{\left(a-\frac 1 2\right)-1}e^{-b\tau_j}
	\Bigg\}\Bigg\{
	\prod\limits_{j=1}^{p} ~\tilde{\tau}_j^{\left(a-\frac 1 2\right)-1}e^{-b\tilde{\tau}_j}
	\Bigg\}\prod\limits_{j=1}^{p}{\beta_j^2\prod\limits_{j=2}^{p}\tilde{\beta}_j^2}~f_3(\bm{\beta},\tilde{\bm{\tau}}) }{\left(f_2(\bm{\beta},\tilde{\bm{\tau}})+1\right)^{n+4p+2\alpha}\prod\limits_{j=2}^{p}\left(1+\frac{\left(\tilde{{\beta}}_j-\mu_j\right)^2}{\nu_j\epsilon_j^2}\right)^{n+4p+2\alpha}}~\nonumber\\
& &~~~~\times~\int_{\mathbb{R}_{+}}\frac{
	\tau_1^{\left(a-\frac 1 2\right)-1}e^{-b\tau_1}
	\left({2\lambda+\frac{1}{\tau_1}+\frac{1}{\tilde{\tau}_1}}\right)^{-\frac 3 2}}{\left(   	\sum\limits_{j=2}^{p}\frac{\beta_j^2}{\tau_j}+\frac{\beta_1^2}{\tau_1}+1\right)^{n+4p+2\alpha}}~d\tau_1\nonumber\\
& & = \infty \nonumber
\end{eqnarray}
for every $\left(\bm{\beta},\tilde{\bm{\tau}},(\tilde{\beta}_2,...,\tilde{\beta}_p)^T,(\tau_2,...,\tau_p)^T\right) \in \mathbb{R}^p \times \mathbb{R}^p_{+} \times \mathbb{R}^{p-1} \times \mathbb{R}^{p-1}_{+} .$ The above integral diverges because 
$\int_{\mathbb{R}_{+}}{
	\tau_1^{(a-\frac 1 2)-1}e^{-b\tau_1}
}d\tau_1 =\infty$ if $a\leq \frac 1 2$ and \[\underset{\tau_1\rightarrow \infty}{\lim}~
\bigfrac{\frac{
		\tau_1^{\left(a-\frac 1 2\right)-1}e^{-b\tau_1}
		\left({2\lambda+\frac{1}{\tau_1}+\frac{1}{\tilde{\tau}_1}}\right)^{-\frac 3 2}}{{\left(   	\sum\limits_{j=2}^{p}\frac{\beta_j^2}{\tau_j}+\frac{\beta_1^2}{\tau_1}+1\right)^{n+4p+2\alpha}}}}{\tau_1^{(a-\frac 1 2)-1}e^{-b\tau_1}
}=\underset{\tau_1\rightarrow \infty}{\lim}\frac{
\left({2\lambda+\frac{1}{\tau_1}+\frac{1}{\tilde{\tau}_1}}\right)^{-\frac 3 2}}{{\left(   	\sum\limits_{j=2}^{p}\frac{\beta_j^2}{\tau_j}+\frac{\beta_1^2}{\tau_1}+1\right)^{n+4p+2\alpha}}}=\frac{
\left({2\lambda+\frac{1}{\tilde{\tau}_1}}\right)^{-\frac 3 2}}{\left(   	\sum\limits_{j=2}^{p}\frac{\beta_j^2}{\tau_j}+1\right)^{n+4p+2\alpha}}\in (0,\infty)
\]

\noindent
\hfill$\Box$

\section{Construction of the sandwich Markov chain} \label{sec:Haarsandwichchain}

\noindent
The two-block Markov chain $\Phi$ can be interpreted as a Data Augmentation (DA) 
algorithm, with $(\bm{\beta}, \sigma^2)$ as the parameter block of interest, and 
$\bm{\tau}$ as the augmented block. The DA algorithm can suffer from slow 
convergence (just like the EM algorithm, its analogous version in likelihood 
maximization). The sandwich algorithm, introduced in 
\cite{Liu:Wu:1999, Hobert:Marchev:2008}, aims to improve the 
convergence and efficiency of the DA algorithm by adding an inexpensive extra 
step in between the two conditional draws of the DA algorithm. In fact, there 
are many DA chains (see \cite{Liu:Wu:1999, Meng:vanDyk:1999, 
Marchev:Hobert:2004, HRR:2011, PKH:2016}, for example) where sandwich 
chains have been constructed and shown to be significantly more efficient with roughly 
the same computational effort per iteration. In this section, we will focus on deriving 
the Haar PX-DA sandwich algorithm in the Bayesian lasso setting. The 
Haar PX-DA algorithm has been shown in \cite{Hobert:Marchev:2008} to be the best 
among a class of sandwich algorithms in terms of efficiency and operator norm. 

A key ingredient in constructing the Haar PX-DA algorithm is a unimodular group 
which acts on the augmented variable space ($\mathbb{R}^p$ in our case). We 
consider the multiplicative group $G$ of positive real numbers, which acts on an 
element of $\mathbb{R}^p$ through scalar multiplication. In particular, if $g \in G$ 
and $\bm{\tau} \in \mathbb{R}^p$, then the result of the action of $g$ on $\bm{\tau}$ 
is given by $g \bm{\tau} = (g \tau_1, g \tau_2, \cdots, g \tau_p)$. Another choice that 
we need to make is the choice of the multiplier function $\chi: G \rightarrow 
\mathbb{R}_+$, which satisfies 
$$
\chi(g_1 g_2) = \chi(g_1) \chi(g_2) 
$$

\noindent
for any pair $g_1, g_2 \in G$, and 
$$
\chi(g) \int_{\mathbb{R}^p} \phi(g \bm{\tau}) d \bm{\tau} = \int_{\mathbb{R}^p} 
\phi(\bm{\tau}) d \bm{\tau} 
$$

\noindent
for any $g \in G$ and any function $\phi: \mathbb{R}^p \rightarrow \mathbb{R}$. In this 
setting, the function $\chi(g) = g^p$ serves as a valid multiplier function. Also, the 
unimodular group $G$ has a Haar measure $\mathcal{H} (dg) = \frac{dg}{g}$. With 
these ingredients in hand, we define the density $f_G$ on $G$ (with respect to the 
Haar measure) by 
\[f_{\mathcal{G}}(g) = \frac{\pi(g\bm{\tau} \mid {\bf Y})\chi(g)}{m(\bm{\tau})}\mathcal{H}
(dg), \]

\noindent
where 
$m(\bm{\tau})= \int_{\mathcal{G}}\pi(g\bm{\tau})\mathcal{\chi}(g)\mathcal{H}(dg)
$ is the normalizing constant. From (\ref{eq:joint}), it follows that 
\[\pi(\bm{\tau} \mid {\bf Y})\propto \frac{\prod_{j=1}^{p}\tau_j^{a-\frac 1 2 - 1 }e^{-b
\tau_j}}{\left\{{\bf{y}}^T{\bf{y}}-{\bf{y}}^TX(X^TX+ D_{\bm{\tau}}^{-1})^{-1}X^T{\bf{y}}
+2\xi\right\}^{\frac n 2 +\alpha}|X^TX+D_{\bm{\tau}}^{-1}|^{\frac 1 2}  }\] 

\noindent
and 
\[f_{\mathcal{G}}(g) \propto 
\frac{g^{p(a-\frac 1 2)-1}e^{-g(\sum_{j=1}^{p}b\tau_j)}}{\left\{{\bf{y}}^T{\bf{y}}-{\bf{y}}^TX(X^TX+\frac 1 g D_{\bm{\tau}}^{-1})^{-1}X^T{\bf{y}}+2\xi\right\}^{\frac n 2 +\alpha}|X^TX+\frac{1}{g}D_{\bm{\tau}}^{-1}|^{{\frac 1 2}}}. 
\]

\noindent
Altohugh $f_{\mathcal{G}}$ is not a standard density, samples from this univariate 
density can be easily generated using a rejection sampling algorithm. Using 
$f_{\mathcal{G}},$ we can now define the Haar PX-DA sandwich Markov chain, 
denoted by $\Phi^{*} = \{(\bm{\beta}_m,\sigma_m^2)\}_{m=0}^{\infty},$ whose one step 
transition from $(\bm{\beta}_m,\sigma_m^2)$ to 
$(\bm{\beta}_{m+1},\sigma_{m+1}^2)$ can be described as follow. \\

1. Draw $\bm{\tau}$ from the distribution $\pi(\cdot~|~\sigma_m^2,\bm{\beta}_m,{\bf 
Y})$ \\

2. Draw $g$ according to the density $f_{\mathcal{G}}.$\\

3. Draw $(\sigma_{m+1}^2,\bm{\beta}_{m+1})$ by the following procedure\\

$~~~~~~~~~~$(a) Draw $\sigma_{m+1}^2$ from $\pi(\cdot ~|~ g\bm{\tau}, {\bf Y}).$\\

$~~~~~~~~~~$(b) Draw $\bm{\beta}_{m+1}$ from $\pi(\cdot ~|~ g\bm{\tau}, \sigma_{m+1}^2,{\bf Y}).$ 

\medskip

\noindent
The lemma below, regarding spectral properties of the Haar PX-DA chain, follows by 
combining Theorem \ref{thm1} and results from \cite{Khare:Hobert:2011}. 
\begin{lemma} \label{lem1}
The operator corresponding to the Haar PX-DA Markov chain $\Phi^*$ is trace class. 
Also, if $\{\lambda^*_i\}_{i=1}^\infty$ and $\{\lambda_i\}_{i=1}^\infty$ are the ordered 
eigenvalues corresponding to $\Phi^*$ and $\Phi$ respectively, then $\lambda^*_i 
\leq \lambda_i$ for $i \geq 1$ with a strict inequality holding for at least one $i$. 
\end{lemma}

\section{Examples} \label{sec:experiments}

\noindent
In this section, we consider two simulated data examples (one each for $n>p $ and 
$n< p $) and a real data example to compare the performance the three block, two 
block and Haar PX-DA sandwich chains. 
 
\subsection{Simulation}

\noindent
We consider a setting with $n=10 < p=15$ for the first simulation, and $n=15 > p=10  $ 
for the second simulation. for both cases,the elements of the design matrix $X$ and 
response $y$ were chosen by generating {i.i.d.} $\mathcal{N}(0,1) $ random variables. 
We fit the Normal-Gamma model in (\ref{eq:model}) with hyper parameters 
$a=0.75,b=2,\xi=100,\alpha=0.$ To compare the efficiency performance of the 
Markov chains, we compute the autocorrelations (up to lag 10) for all the Markov 
chains for the function $\left({\bf Y}-X\bm{\beta}\right)^T\left({\bf Y}-X\bm{\beta}\right)+
\sigma^2 .$ The results are summarized in Table 1 and Figure 1 for the first simulation, 
and in Table 2 and Figure 2 for the second simulation. We can clearly see that for both 
datasets, the two block Gibbs sampler has significantly lower autocorrelations than the 
three block Gibbs sampler, and that the magnitude of the autocorrelations for the 
sandwich Markov chain is lowest. 
	\begin{table}[H]
		\caption{First ten autocorrelations for simulated data with $n < p$ } 
		\centering\begin{tabular}{rrrrrrrrrrr}
			\toprule
			\bfseries Lag & 1 & 2 & 3 & 4 & 5 & 6 & 7 & 8 & 9& 10 \\\midrule
		\bfseries Three block
			& 0.5018& 0.2687& 0.1510& 0.0794& 0.0403& 0.0255& 0.0104& 0.0053& 0.0050& -0.0101
		 \\\midrule
		\bfseries Two Block
		& 0.1645& 0.0725& 0.0550& 0.0223& 0.0048& 0.0093& 0.0158& 0.0034& 0.0166& 0.0174
        \\\midrule
		\bfseries Sandwich
		& 0.0177& 0.0010& -0.0078& -0.0362& -0.0111& -0.0152& -0.0068& 0.0108& 0.0166 & 0.0053
		\\\bottomrule
		\end{tabular}
	\end{table}

\begin{figure}[H]
	\centering
	\includegraphics[height=8cm]{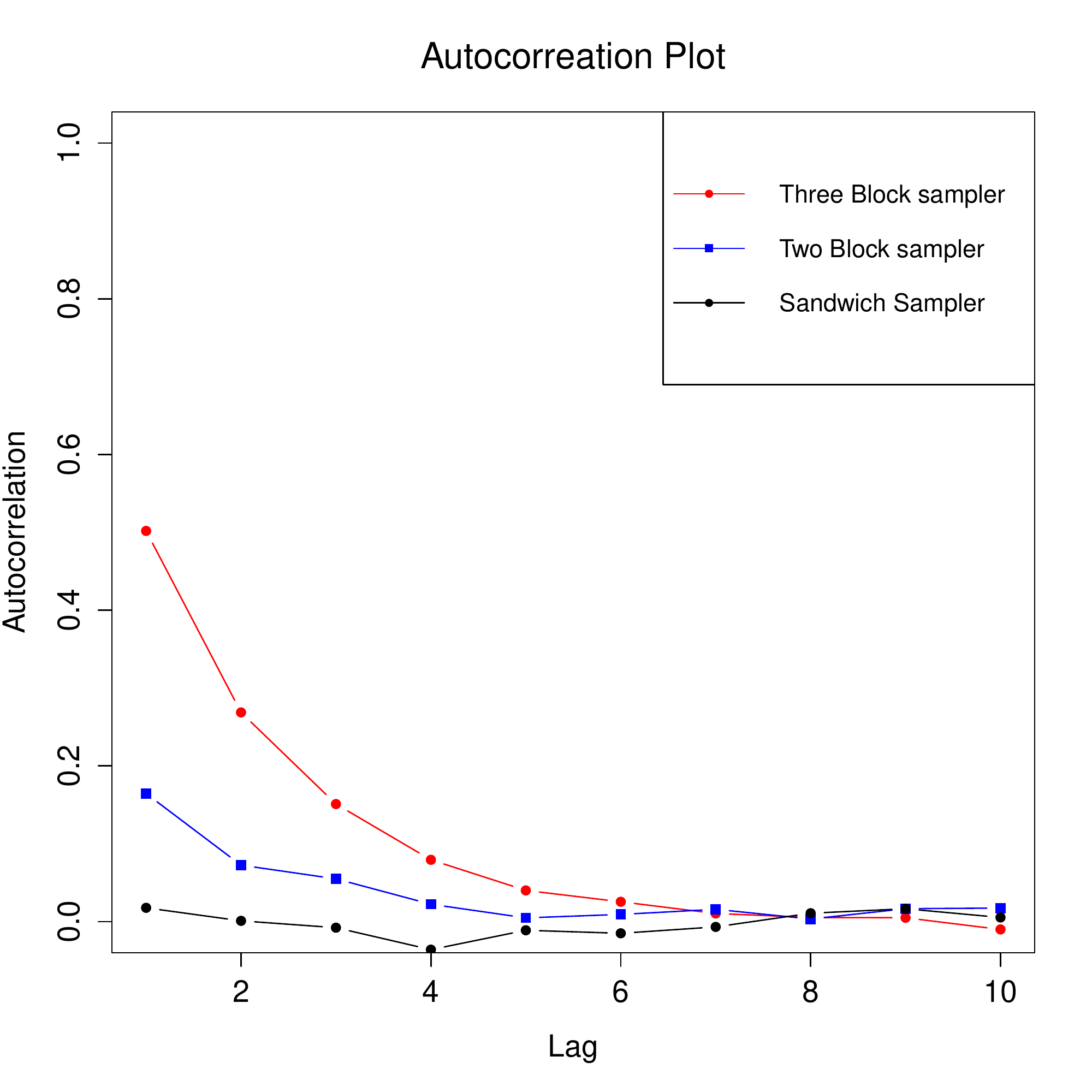}
	\caption{Autocorrelation plot of $\left({\bf Y}-X\bm{\beta}\right)^T\left({\bf Y}-X\bm{\beta}\right)+\sigma^2$ for simulated data with $n<p.$ }
\end{figure}
	
\begin{table}[H]
			\caption{First ten autocorrelations for simulated data with $n > p$ } 
			\centering\begin{tabular}{rrrrrrrrrrr}
				\toprule
				\bfseries Lag & 1 & 2 & 3 & 4 & 5 & 6 & 7 & 8 & 9& 10 \\\midrule
				\bfseries Three block
				& 0.5413& 0.3497& 0.2308& 0.1562& 0.1033& 0.0638& 0.0423& 0.0231& 0.0197& 0.0143 \\\midrule
				\bfseries Two Block
				& 0.1155&  0.0544&  0.0297&  0.0413&  0.0340&  0.0130&  0.0043&  -0.0022&  0.0182&  -0.0045
				\\\midrule
				\bfseries Sandwich
				&-0.0198& -0.0238& -0.0077& 0.0512& -0.0079& 0.0232& -0.0036& -0.0005& -0.0048& -0.0259
				\\\bottomrule
			\end{tabular}
\end{table}
	
	\begin{figure}[H]
		\centering
		\includegraphics[height=8cm]{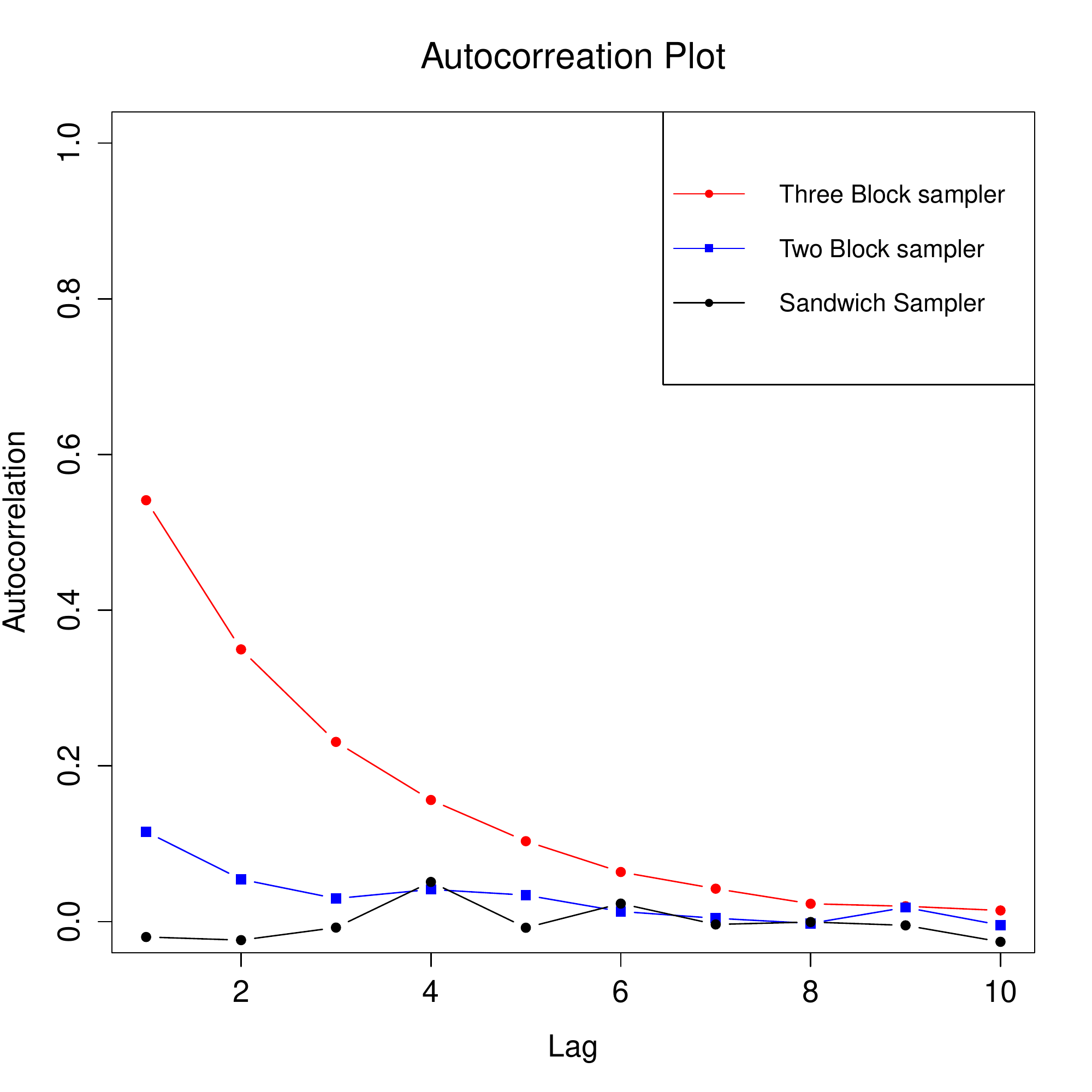}
		\caption{Autocorrelation plot of $\left({\bf Y}-X\bm{\beta}\right)^T\left({\bf Y}-X\bm{\beta}\right)+\sigma^2$ for simulated data with $n>p.$ }
	\end{figure}

\subsection{Real data example}

\noindent
We consider the wheat data set from 
Perez and de los Campos \cite{Perez:delosCampos:2014}, which is available in 
the R package BGLR. The data was obtained from numerous international trials 
across a wide variety of wheat-producing environments. For our analysis, we consider 
the average grain yield for a particular environmental condition (there are four to 
choose from) as the response variable, and 20 binary variables containing genotypic 
information as the predictors. We  fit the Normal-Gamma model in (\ref{eq:model}) 
with $a=0.75,b=0.2,\xi=1,\alpha=0$ and compute autocorrelations for the 
function $\left({\bf Y}-X\bm{\beta}\right)^T\left({\bf Y}-X\bm{\beta}\right)+\sigma^2$ for 
the three block, two block and Haar PX-DA sandwich chains. The results are shown in 
Table 3 and Figure 3. As in the simulated data examples, the two-block chain has 
lower autocorrelations than the three-block chain, and the Haar PX-DA sandwich chain 
is the most efficient among all three Markov chains. 
	\begin{table}[H]
		\caption{First ten autocorrelations with wheat data } 
		\centering\begin{tabular}{rrrrrrrrrrr}
			\toprule
			\bfseries Lag & 1 & 2 & 3 & 4 & 5 & 6 & 7 & 8 & 9& 10 \\ \midrule
		\bfseries Three block
	&	0.4579 & 0.2246& 0.1193& 0.0433& 0.0214& 0.0194& 0.0166 & 0.0114& -0.0072& -0.0068\\\midrule
	\bfseries Two block
		& 0.0799& 0.0159& 0.0053& 0.0065& 0.0072& 0.0015& 0.0022& -0.0224& -0.0103& -0.0131\\\midrule
	\bfseries Sandwich
		& 0.0535 & 0.0163& 0.0009 & 0.0344& 0.0073& 0.0308& -0.0190& -0.0438& -0.0098& -0.0161 \\\bottomrule
		\end{tabular}
	\end{table}
	
\begin{figure}[H]
	\centering
	\includegraphics[height=8cm]{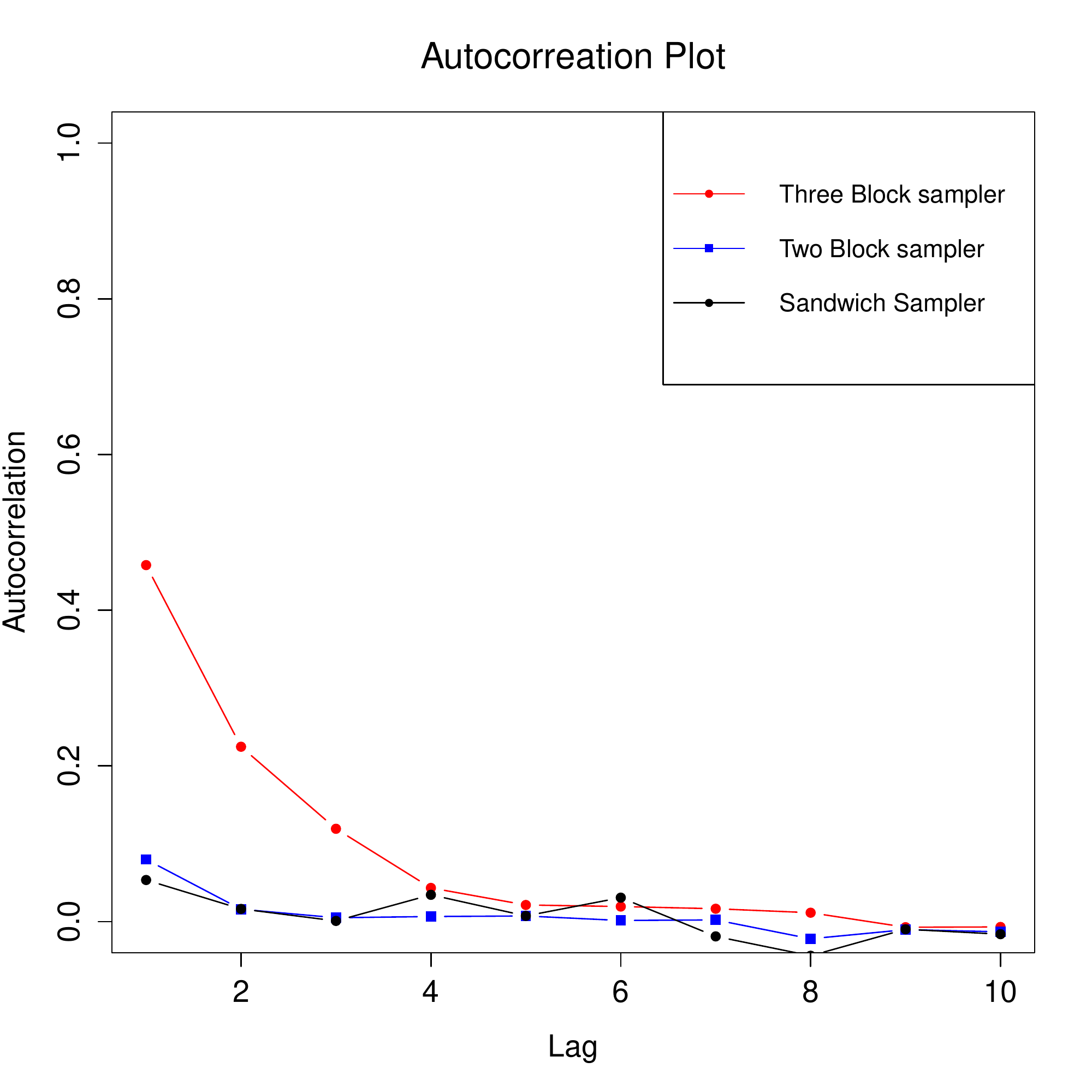}
	\caption{Autocorrelation plot of $\left({\bf Y}-X\bm{\beta}\right)^T\left({\bf Y}-X\bm{\beta}\right)+\sigma^2$ for the wheat data. }
\end{figure}

\subsection{Discussion of numerical results}

\noindent
For both the simulated and real data settings, the two-block chain clearly has a 
significantly better performance than the three-block chain. For example, in all the 
settings the Lag 1 autocorrelation drops by 80\% or more when we compare the 
three-block and the two-block chains. These findings support the theoretical results 
(Theorem \ref{thm1} and Theorem \ref{thm2}) in the paper. Since the two-block chain 
and the three-block chain require the same computational effort, our theoretical and 
experimental results, support the overall conclusion that a practitioner should prefer 
the two-block chain over the three-block chain. 

The comparison between the two-block and the Haar PX-DA chain is not as decisive. 
We elaborate on this below. Clearly, in all the experimental settings, the sandwich 
algorithm performs better than the two-block chain, and substantially so for the 
simulated datasets. These findings support the theoretical results (Lemma \ref{lem1}). 
However, note that the Haar PX-DA chain requires an extra computational step 
(sampling from $f_{\mathcal{G}}$) as compared to the two-block chain. This extra 
step requires a univariate rejection sampler. We have found that the performance 
of this rejection sampler varies with the choice of $a, b, \xi$ and $\alpha$. In the 
best case, it takes twice as much time as the other steps, and in the worst case, can 
take more than 10 times as much time as the other steps. Of course, the performance 
gains can sometimes easily offset this computational overhead, but a practitioner 
should make a careful determination in their specific setting.

\section*{Appendix}

\newtheorem{proposition}{Proposition A \hspace{-0.1in}}
\begin{proposition}
Let $x \sim N(\mu,\sigma^2) ,$ then $~\int\limits_{0}^{c} \frac 1 x e^{-\frac{(x-\mu)^2}
{2\sigma^2}}\,dx =\infty~$ for any positive constant $c>0$. 
\end{proposition}

\noindent
{\it Proof} The result follows by noting that 
$$
x\cdot \abs*{\frac 1 x e^{-\frac{(x-\mu)^2}{2\sigma^2}}} \stackrel{~x\rightarrow 0^{+}}
{\longrightarrow} e^{-\frac{\mu^2}{2\sigma^2}} > 0, 
$$ 

\noindent
and $\int\limits_{0}^{c} \frac{1}{x} \,dx = \infty$. 
\hfill$\Box$ 

\begin{proposition}
Let $x \sim N(\mu,\sigma^2) ,$ then $~\int\limits_{0}^{c_1} \frac 1 {x(c_2-\ln x)}~ e^{-
\frac{(x-\mu)^2}{2\sigma^2}}\,dx =\infty~$ for any positive constant $c_1$ and 
$c_2 > \ln c_1$. 
\end{proposition}

\noindent
{\it Proof} We assume without loss of generality that $c_1 < 1$. The result follows by 
noting that 
$$ 
\frac{\frac 1 {x(c_2-\ln x)}~ e^{-\frac{(x-\mu)^2}{2\sigma^2}}}{-\frac{1}{x\ln x}}
=\frac{-\ln x}{c_2-\ln x}~ e^{-\frac{(x-\mu)^2}{2\sigma^2}}\stackrel{~x\rightarrow 0^{+}}
{\longrightarrow} e^{-\frac{\mu^2}{2\sigma^2}} \in (0,\infty), 
$$

\noindent
and 
$$
~\int\limits_{0}^{c_1}-\frac 1 { x \ln x} \,dx = \int\limits_{0}^{c_1}-\frac 1 {\ln x}\,d\ln x~ 
\stackrel{t=\ln x}{=\joinrel=\joinrel=\joinrel=} \int\limits_{-\infty}^{\ln c_1} -\frac 1 t \,dt =
\int\limits_{-\ln c_1}^{\infty} \frac 1 t \,dt = \infty. 
$$

\noindent
\hfill$\Box$

\begin{proposition}
Suppose the trandom variable U has a t-distribution with scale parameter $\kappa$ , 
location parameter $\vartheta$ and degrees of freedom $\nu.$ Then for $\nu>2,$
\[E(U^2) \geq \frac{k^2\nu}{\nu-2}\]. 
\end{proposition}

\noindent 
{\it Proof} If the random variable $U$ has a t-distribution with scale parameter $
\kappa,$ location parameter $\vartheta$ and degrees of freedom $\nu$ then $U=
\vartheta+\kappa T$ where $T$ is a standard t-distribution with $\nu$ degrees of 
freedom. Hence
\begin{align*}
E(U^2)=E(\vartheta+\kappa T)^2&=E(\vartheta^2+\kappa^2T^2+2\vartheta\kappa T)\\
&\geq E(\kappa^2T^2+2\vartheta\kappa T)\\
&= \kappa^2E(T^2)+2\vartheta \kappa E( T)\\
& =\kappa^2(Var(T)+E^2(T))\\
&= \frac{\kappa^2\nu}{\nu-2}. 
\end{align*}

\noindent
\hfill$\Box$

\begin{proposition}
Let $\mu_1=e_1^T\mu,$ i.e. the first component of $\mu$, $\nu_1=2n+8p+4\alpha-1$, 
and $\epsilon_1= \sqrt{\frac{f_2(\bm{\beta},\tilde{\bm{\tau}})+1}{(2\lambda+\frac{1}
{\tau_1}+\frac{1}{\tilde{\tau}_1})\nu_1}}$. Then there is a finite constant 
$f_3(\bm{\beta},\tilde{\bm{\tau}}) $ such that 
$$
\int_{\mathbb{R}}\frac{\tilde{\beta}_1^2} {\left(1+\frac{\left(\tilde{\bm{\beta}}_1-
\mu_1\right)^2}{\nu_1\epsilon_1^2}\right)^{\frac{1+\nu_1}{2}}}d{\tilde{{\beta}}_1} \geq 
f_3(\bm{\beta},\tilde{\bm{\tau}}) \left(\frac{1}{2\lambda+\frac{1}{\tau_1}+\frac{1}
{\tilde{\tau}_1}}\right)^{\frac 3 2}. 
$$
\end{proposition}

\noindent
{\it Proof} Note that
\[\int_{\mathbb{R}}\frac{\tilde{\beta}_1^2} {\left(1+\frac{\left(\tilde{\bm{\beta}}_1-
\mu_1\right)^2}{\nu_1\epsilon_1^2}\right)^{\frac{1+\nu_1}{2}}}d{\tilde{{\beta}}_1}=
\epsilon_1\frac{\Gamma(\frac{\nu_1}{2})\sqrt{\pi\nu}}{\Gamma(\frac{\nu+1}{2})}E(U^2),
\]
where $U$ follows a t-distribution with scale $\epsilon_1,$ location $\mu_1$ and 
degrees of freedom $\nu_1.$ Using Proposition A3, we get that

\begin{align*}
\int_{\mathbb{R}}\frac{\tilde{\beta}_1^2} {\left(1+\frac{\left(\tilde{\bm{\beta}}_1-
\mu_1\right)^2}{\nu_1\epsilon_1^2}\right)^{\frac{1+\nu_1}{2}}}d{\tilde{{\beta}}_1} & =
\epsilon_1\frac{\Gamma(\frac{\nu_1}{2})\sqrt{\pi\nu}}{\Gamma(\frac{\nu+1}{2})}E(U^2)\\ 
& \geq \frac{\Gamma(\frac{\nu_1}{2})\sqrt{\pi\nu}}{\Gamma(\frac{\nu+1}{2})}\frac{\nu_1}
{\nu_1-2}~\epsilon_1^3,\\
& = \frac{\Gamma(\frac{\nu_1}{2})\sqrt{\pi\nu}}{\Gamma(\frac{\nu+1}{2})}\frac{\nu_1}
{\nu_1-2}\left(\frac{f_2(\bm{\beta},\tilde{\bm{\tau}})+1}{\nu_1}\right)^{\frac 3 2} 
\left(\frac{1}{2\lambda+\frac{1}{\tau_1}+\frac{1}{\tilde{\tau}_1}}\right)^{\frac 3 2}\\
& =f_3(\bm{\beta},\tilde{\bm{\tau}}) 
\left({2\lambda+\frac{1}{\tau_1}+\frac{1}{\tilde{\tau}_1}}\right)^{-\frac 3 2}.
\end{align*}

\noindent
where $f_3(\bm{\beta},\tilde{\bm{\tau}})=\frac{\Gamma(\frac{\nu_1}{2})\sqrt{\pi\nu}}
{\Gamma(\frac{\nu+1}{2})}\frac{\nu_1}{\nu_1-2}\left(\frac{f_2(\bm{\beta},
\tilde{\bm{\tau}})+1}{\nu_1}\right)^{\frac 3 2}$. \hfill$\Box$

\end{document}